\documentclass[11pt]{article}

\usepackage{mathtools}
\usepackage{amssymb}
\usepackage{amsthm}
\usepackage{xcolor}
\usepackage{tikz}
\usepackage{mathrsfs}
\usepackage{marginnote}

\usepackage{geometry}
\definecolor{myteal}{HTML}{004b6f}
\definecolor{citegrey}{HTML}{888888}
\usepackage{hyperref}
\hypersetup{
    colorlinks=true,
    linkcolor=myteal,
    citecolor=citegrey,
    urlcolor=myteal,
}

\newcommand{\mysetminusD}{\hbox{\tikz{\draw[line width=0.6pt,line cap=round] (3pt,0) -- (0,6pt);}}}
\newcommand{\mysetminusT}{\mysetminusD}
\newcommand{\mysetminusS}{\hbox{\tikz{\draw[line width=0.45pt,line cap=round] (2pt,0) -- (0,4pt);}}}
\newcommand{\mysetminusSS}{\hbox{\tikz{\draw[line width=0.4pt,line cap=round] (1.5pt,0) -- (0,3pt);}}}

\renewcommand{\setminus}{\mathbin{\mathchoice{\mysetminusD}{\mysetminusT}{\mysetminusS}{\mysetminusSS}}}

\DeclarePairedDelimiter{\abs}{\lvert}{\rvert}
\DeclarePairedDelimiter{\Angle}{\langle}{\rangle}
\renewcommand{\angle}{\Angle}
\DeclarePairedDelimiter{\norm}{\lVert}{\rVert}

\newcommand{\R}{\mathbb{R}}
\newcommand{\N}{\mathbb{N}}
\newcommand{\E}{\mathbb{E}}
\renewcommand{\P}{\mathbb{P}}

\renewcommand{\Pr}{\mathrm{Pr}}

\newcommand{\diff}{\mathrm{d}}
\newcommand{\indic}{\mathbf{1}}
\DeclareMathOperator{\Card}{Card}

\newtheorem{theorem}{Theorem}
\newtheorem{lemma}{Lemma}
\newtheorem{proposition}{Proposition}
\newtheorem{corollary}{Corollary}

\theoremstyle{definition}
\newtheorem{assumption}{Assumption}

\newtheorem{example}{Example}

\usepackage{authblk}

\newlength{\affilskip}
\author[1]{Félix Foutel-Rodier}
\author[2]{Emmanuel Schertzer}

\affil[1]{
    Université Paris Cité, CNRS, MAP5, F-75006 Paris, France
    \vspace{\affilskip}
}
\affil[2]{
    Faculty of Mathematics, University of Vienna,\authorcr Oskar-Morgenstern-Platz 1, 1090 Wien, Austria
    \vspace{\affilskip}
}

\title{Non-linear branching processes and Crump--Mode--Jagers processes
with interaction}

\begin{document}

\maketitle

\begin{abstract}
    We consider a class of Crump--Mode--Jagers processes with
    interaction, constructed by removing a newly born offspring
    with a probability that depends on the age structure of the
    population at its birth time. We prove a law of large numbers for the
    tree structure of the process in a local topology, and show how this
    result condenses several other limit theorems (convergence of the
    empirical age distribution, of ancestral lineages). Beyond this
    specific example, our work illustrates a more general principle that
    we formalise.
    As in standard propagation of chaos, the trees generated by typical
    individuals become independent as the number of individuals goes to
    infinity. This allows us to express the distribution of the local
    tree structure around a typical individual in terms of a
    time-inhomogeneous branching process, which we call a non-linear
    branching process.
\end{abstract}

\section{Introduction}

\subsection{Motivation}

Crump--Mode--Jagers (CMJ) processes are branching particle systems in
continuous time in which the ages at which individuals produce offspring
are distributed as a general point process \cite{jagers1975branching,
Nerman1981, jagers1984growth}. This great generality provides a
very flexible framework to model age-structured populations,
which has found recent applications to epidemiology \cite{britton2019estimation, 
favero2022modelling, foutel2022individual, pakkanen2023unifying}. However, 
in this context, the assumption that individuals reproduce independently
becomes limiting when modeling several features of real-world outbreaks.
A first simple example is that, if the disease provides immunity after
infection, the pool of individuals susceptible to the disease decreases
during the spread of the epidemic, leading to a reduction in the
transmission rate. This effect generates the typical wave shape observed
in SIR-type models and cannot be modeled appropriately by a branching process.
Even in the early phase of the epidemic, when this effect is arguably
negligible, the branching property might fall short. The onset of an
epidemic often leads to behavioural adaptations and to the enforcement of
control measures (lockdown, instructions given to the general population)
which influence the transmission rate. This creates a dependence between
infections through the adaptation of the population to the epidemic which
disrupts the branching property. These examples lead us to introduce a
class of CMJ processes with general interaction, that can encompass both
situations. Our main result is a law of large numbers for the tree
structure generated by the births in these CMJ processes. 

Although the results we derive on CMJ processes are interesting for
themselves, the main aim of our work is to illustrate a broader approach
to studying population processes with interactions. In many population
models, individuals only interact through empirical averages over
the population, sometimes referred to as \emph{mean-field interaction}.
For instance, individuals interact through local population density in
spatial ecology, allelic frequencies in population genetics, or density
of infected individuals in epidemiology. A very similar situation is
encountered in many physical particle systems, for which a rich theory
has been developed around the notion of \emph{propagation of chaos}.
Briefly, if such a system starts from an initial configuration where
particles are independent, it is expected that particles remain
asymptotically independent as their number goes to infinity. As a
consequence, any macroscopic state of the system follows a law of large
numbers and converges to a deterministic limit which can be expressed in
terms of the law of one typical particle. This typical particle is
distributed as a so-called \emph{non-linear (Markov) process}, whose law
can be obtained by solving a non-linear equation. We refer the reader to
\cite{sznitman1991, Chaintron2022} and references therein for a detailed
exposition of these notions and a historical perspective.

One specificity of population models compared to most physical systems to
which propagation of chaos applies is that mass (that is the number of
particles) is not conserved: individuals reproduce and die. This feature
of biological systems might be seen as a technical obstacle but, more
interestingly, it creates a branching structure which is absent from the
standard theory of mean-field particle systems. The main conceptual
contribution of our work is to show that propagation of chaos applies to
this branching structure, and that a law of large numbers should hold
generally in the sense of local weak convergence of graphs. More
precisely, we expect the trees grown out of two typical initial particles
to become asymptotically independent and each distributed as a
time-inhomogeneous branching process. Following the terminology of
propagation of chaos, we call this limiting tree a \emph{non-linear
branching process} as its distribution solves a non-linear equation. It
can be used to describe the limit of the local tree structure around a
typical individual in the population. This connection enables the rich
literature on branching processes to be used to study the branching
structure of large interacting populations, in particular the successful
spinal decomposition techniques \cite{lyons1995conceptual, Shi2015,
marguet2019uniform}. We will illustrate these principles on our CMJ
example, but we expect that they apply broadly.

The idea of approximating a Markov branching process with mean-field
interaction by a time-inhomogeneous one appears in several
works \cite{dewitt2024mean, overbeck1996nonlinear, henry2023time,
calvez2022dynamics}. In particular, \cite{overbeck1996nonlinear}
defines a general notion of non-linear superprocess, which can be
seen as the measure-valued version of our non-linear branching
process. Comparison to a branching process is also the key step in
\cite{calvez2022dynamics, henry2023time} to study the ancestral
lineages in a population with interaction.

The contribution of our work in this respect is that we lift the
comparison with a time-inhomogeneous process to tree structures.
Beyond allowing us to study the genealogy of the population over
short times, this extension is the crucial step that enables us to
add interaction to CMJ processes. These CMJ processes have received a
lot of attention in the past after the seminal work of Jagers and
Nerman \cite{Nerman1981, Nerman1984, jagers1984growth}. The absence of
the standard Markov property for those objects requires new techniques
compared to existing methods for populations with interaction, which are
based on Poisson-driven stochastic differential equations. In comparison
to previous work that heavily rely on stochastic calculus, our
approach is based on a coupling technique tailored for CMJ processes.
See the discussion after Theorem~\ref{thm:main} for more details, and
Section~\ref{SS:fullCoupling} for an explicit description of our coupling
technique.

\subsection{Crump--Mode--Jagers processes with interaction}
\label{SS:constructionCMJ}

We start by introducing CMJ processes without interaction, and then show
how to incorporate interaction by removing some individuals with a
probability that depends on the state of the population. The construction
of a CMJ process is based on a random point measure $\mathcal{P}$ on
$(0,\infty)$ whose atoms encode the ages at which a typical individual
gives birth. For simplicity, we assume that the atoms of $\mathcal{P}$ are
almost surely distinct and let $\abs{\mathcal{P}}$ denote the mass of
$\mathcal{P}$. Following the Ulam--Harris convention, we will label individuals
in the population as finite words, which are elements of the set
\[
    \mathcal{U} = \{ \varnothing \} \cup \bigcup_{n \ge 1} \N^n.
\]
Our process will start from a large number of individuals and is a forest
rather than a single tree. In that direction, let 
\[
    \forall i \ge 1,\quad \mathcal{U}_i = \{ (i,u) : u \in \mathcal{U} \}
\]
be the set of individuals descended from $i$, which we interpret as the
$i$-th ancestor in the population. We use the standard notation $\abs{u}$
for the length of the word $u$ (interpreted as the generation of $u$) and
write $u \preceq v$ if $u$ is an ancestor of $v$. 

Let $(\mathcal{P}_u)_{u \in \mathcal{U}}$ be an i.i.d.\ collection of
point processes distributed as $\mathcal{P}$ and fix some probability 
density $g$ on $(0,\infty)$, corresponding to the distribution of initial
ages in the population. Let $(A^0_i)_{i \ge 1}$ be a sequence of i.i.d.\
random variables distributed as $g(a)\diff a$ and define 
\begin{equation} \label{eq:initialCondition}
    \forall i \ge 1,\quad \sigma_i = -A^0_i,
    \qquad 
    \mathcal{P}^0_i = \sum_{a \in \mathcal{P}_i} \indic_{\{ a \ge A^0_i\}}
    \delta_a.
\end{equation}
The variable $\sigma_i$ is the birth time (in the past) of individual
$i$, and $\mathcal{P}^0_i$ is obtained by shifting $\mathcal{P}_i$ by
$\sigma_i$ and removing all atoms at negative times. For every $u\in{\cal U}$,
let us denote by $A_{u,1} < A_{u,2} < \dots$ the atoms of $\mathcal{P}_u$
(or of $\mathcal{P}^0_i$ if $u=i$) in increasing order, with the convention
that $A_{u,k} = \infty$ if $k > \abs{\mathcal{P}_u}$. Inductively, we let
\[
    \forall u \in \mathcal{U} \setminus \{ \varnothing \},\quad 
    \sigma_{ui} = \sigma_u + A_{u,i}.
\]
In words, each individual reproduces at ages given by the atoms of an
independent copy of $\mathcal{P}$, and its offspring is labeled in
increasing order of birth times. For each $i \ge 1$, this procedure
constructs a collection of birth times $\mathcal{T}_i \coloneqq
(\sigma_u)_{u \in \mathcal{U}_i}$ that we call a \emph{CMJ tree}. 

Interaction is introduced in the population by removing an individual
with a probability that depends on the state of the population at its
birth time. 
Fix a number $N$ of initial individuals and a function $C \colon \R_+
\times \mathcal{M}(\R_+) \to [0,1]$, where $\mathcal{M}(\R_+)$ is the
space of finite measures on $\R_+$. We construct a pruning of the CMJ
tree such that an individual born at time $t$ is only kept in the
population with probability $C(t, \mu^N_{t^-})$, where $\mu^N_t$ is the
empirical measure of birth times at time $t$.

More formally, let $(\omega_u)_{u \in \mathcal{U}}$ be an i.i.d.\
collection of uniform random variables on $[0,1]$. Define
\[
    \forall u \in \mathcal{U},\qquad 
    \chi^N_u = \indic\big(\{ \omega_u \le C(\sigma_u, \mu^N_{\sigma_u^-})\}\big),
    \qquad
    I^N_u = \prod_{v \preceq u} \chi^N_u,
\]
where $\mu_t^N$ is the empirical distribution of ages in the population:
\[
    \forall t \ge 0,\quad 
    \mu^N_{t}(\diff a) = \frac{1}{N} \sum_{i=1}^N
    \indic_{\{I^N_u=1, \sigma_u \le t\}} \delta_{t-\sigma_u}(\diff a).
\]
The variable $\chi^N_u$ indicates whether the birth of individual $u$ is
`effective' or not, and $I^N_u$ indicates if $u$ or one of its ancestors
has been removed from the population. The population is made of those
individuals such that $I^N_u = 1$. The existence of such a collection
of random variables can be obtained inductively, by considering
individuals in increasing order of their birth times $(\sigma_u)_{u \in
\mathcal{U}}$. Finally, we define
\[
    \forall u \in \mathcal{U},\quad 
    \sigma^N_u =
    \begin{dcases}
        \sigma_u &\text{if $I^N_u = 1$,} \\
        \infty &\text{if $I^N_u = 0$,}
    \end{dcases}
\]
and let $\mathcal{T}^N_i = ( \sigma^N_u )_{u \in \mathcal{U}_i}$ be the
tree rooted at the $i$-th ancestor, after pruning. We call
$\mathcal{T}^N_i$ the $i$-th \emph{CMJ tree with interaction}.

\subsection{Non-linear branching processes}

As alluded to above, our proof adapts ideas from propagation of chaos to
branching particle systems. We will show that, as $N \to \infty$, the
subtrees grown out of different ancestors are asymptotically independent
and distributed as a time-inhomogeneous CMJ tree that we now construct.
Some informal motivation for this definition is provided at the end of
the section.

Fix some deterministic function $(m_t)_{t \ge 0}$, where $m_t \in \mathcal{M}(\R_+)$.
Let $\chi'_u$ and $I'_u$ be constructed as above, but replacing the
empirical distribution of ages $\mu^N_t$ by $m_t$. Namely,
\[
    \forall u \in \mathcal{U},\qquad 
    \chi'_u = \indic\big\{ \omega_u \le C(\sigma_u^-, m_{\sigma_u^-}) \}\big), \qquad
    I'_u = \prod_{v \preceq u} \chi'_u,
\]
and, as above, let
\[
    \forall u \in \mathcal{U},\quad 
    \sigma'_u = 
    \begin{dcases}
        \sigma_u &\text{if $I'_u = 1$,} \\
        \infty &\text{if $I'_u = 0$.}
    \end{dcases}
\]
Let the empirical distribution of birth times of the progeny of a single
ancestor be
\[
    \forall t \ge 0,\quad 
    \mu'_t = \sum_{u \in \mathcal{U}_1} \indic_{\{ \sigma'_u \le t \}} \delta_{t-\sigma_u}.
\]
If the function $(m_t)_{t \ge 0}$ satisfies that 
\begin{equation} \label{eq:McKeanVlasov}
    \forall t \ge 0,\quad
    m_t(\diff a) = \E[ \mu'_t(\diff a)],
\end{equation}
we call the random tree $\mathcal{T}' = (\sigma'_u)_{u \in
\mathcal{U}_1}$ a \emph{non-linear CMJ tree}.

The existence of a non-linear CMJ tree is no longer guaranteed by a
simple induction, but the following result shows that it can be connected
to the existence of a solution to a non-linear renewal equation. An
appropriate notion of weak solution for this equation has been developed
in the monograph \cite{webb1985theory}, which we recall in
Section~\ref{S:PDE} along with a proof of the following result. We use
the convention that, if $f \in \mathrm{L}^1(\R_+)$, we identify it with
the measure $f(a)\diff a$ and simply write $C(t,f)$.

\begin{proposition} \label{prop:nonLinearExistence}
    There exists a non-linear CMJ tree if and only if there is a weak solution
    $(u_t(a);\, t,a \ge 0)$ to 
    \begin{align}
    \begin{split} \label{eq:main}
        (\partial_t + \partial_a) u_t(a) &= 0 \\
        u_t(0) &= C\big(t, u_t\big) \int_0^\infty u_t(a) \tau(a) \diff a\\
        u_0(a) &= g(a).
    \end{split}
    \end{align}
    The empirical distribution of ages of the corresponding non-linear
    tree satisfies
    \[
        \forall t \ge 0,\quad \E\big[ \mu'_t(\diff a) \big] = u_t(a) \diff a.
    \]
\end{proposition}

Let us end this section by providing some intuition about the definition
of a non-linear CMJ tree. We follow closely some seminal ideas in
standard propagation of chaos.
Suppose that a law of large numbers holds and that the empirical measure 
of ages $\mu^N_t(\diff a)$ converges to a deterministic limit $u_t(a)
\diff a$. By construction, provided that $C$ is continuous, a birth at
time $t$ is pruned with a probability that converges to $1-C(t, u_t)$.
Therefore, one can expect that the tree $\mathcal{T}^N_i =
(\sigma^N_u)_{u \in \mathcal{U}_i}$ converges to some limit
$\mathcal{T}^\infty_i = (\sigma^\infty_u)_{u \in \mathcal{U}}$, for which
a birth at time $t$ is kept with probability $C(t, u_t)$. 

Furthermore, the exchangeability of the system entails that
\[
    \E[\mu^N_t(\diff a)] 
    = \frac{1}{N} \sum_{i=1}^N 
      \E\Big[ \sum_{u \in \mathcal{U}_i} \indic_{\{\sigma^N_u \le t\}}
        \delta_{t-\sigma^N_u}(\diff a) \Big]
    = \E\Big[ \sum_{u \in \mathcal{U}_1} \indic_{\{\sigma^N_u \le t\}} 
        \delta_{t-\sigma^N_u}(\diff a) \Big].
\]
Hence, if both sides of this identity converge, we expect that 
$u_t(a)$ and $\mathcal{T}^\infty_1$ are related through
\[
    u_t(a) \diff a 
    = \E\Big[ 
        \sum_{u \in \mathcal{U}_1} 
        \indic_{\{ \sigma^\infty_u \le t\}} \delta_{t-\sigma^\infty_u}(\diff a)
    \Big],
\]
which is exactly our definition of a non-linear CMJ tree, with
$m_t = u_t$.

\subsection{Main result}
\label{sect:main-results}

We will derive our result under the following technical conditions.

\begin{assumption} \label{ass:PDE}
    Let $\tau(\diff a) = \E[\mathcal{P}(\diff a)]$ be the intensity
    measure of the birth point process. We assume that:
    \begin{enumerate}
    \item[(i)] The measure $\tau$ has a density $\tau(\diff a) = \tau(a)\diff a$
        which satisfies that $\tau \in \mathrm{L}^\infty(\R_+)$.
    \item[(ii)] The map $(t, \mu) \mapsto C(t,\mu)$ is continuous (for
        the weak topology) and, for any $T > 0$, there exists $L \equiv
        L(T)$ such that 
        \begin{equation} \label{eq:Lipschitz}
            \forall \mu,\nu \in \mathcal{M}(\R_+), \:
            \forall t \le T,
            \quad \abs{C(t, \mu) - C(t, \nu)} \le L d_{\Pr}(\mu, \nu),
        \end{equation}
        where $d_{\Pr}$ is Prohorov's distance on $\mathcal{M}(\R_+)$.
    \end{enumerate}
\end{assumption}

Recall from the previous section that ${\cal T}_i^N$ is the $i$-th tree
in the CMJ process with interaction. Formally, we encode trees as
elements of 
\begin{equation} \label{eq:defTrees}
    \mathscr{T} \coloneqq 
    \big\{ (s_u)_{u \in \mathcal{U}} \in \big(\R_+ \cup \{\infty\}\big)^{\mathcal{U}}
        : u \preceq v \implies s_u \le s_v    
    \big\}.
\end{equation}
The coordinate $s_u$ is the birth time of individual $u$, with the
interpretation that $s_u = \infty$ if $u$ is not in the tree. We endow
each coordinate with the usual topology on $\R_+$ and treat $\infty$ as
an isolated point. The space $\mathscr{T}$ is further endowed with the
product topology. In a tree interpretation, convergence under this
topology is equivalent to local convergence (around the root) of the
discrete tree structure and of the birth times.
In order to view $\mathcal{T}^N_i$ as a random element of $\mathscr{T}$,
one needs to shift all labels in such a way that the root of the
$\mathcal{T}^N_i$ becomes $\varnothing$ and not $i$. We will always
implicitly do so.

We are now ready to state our main result. Recall that we have
constructed trees as random elements of $\mathscr{T}$ defined in
\eqref{eq:defTrees}, endowed with a local topology.

\begin{theorem} \label{thm:main}
    Suppose that Assumption~\ref{ass:PDE} is satisfied. Then
    \eqref{eq:main} has a unique solution. Let $\mathcal{T}'$ be the
    corresponding non-linear CMJ tree and $\mathscr{L}(\mathcal{T}')$ be
    its distribution. Then
    \[
        \lim_{N \to \infty} \frac{1}{N} \sum_{i=1}^N \delta_{\mathcal{T}_i^N}
        =
        \mathscr{L}(\mathcal{T}'),
    \]
    in law as random measures on $\mathscr{T}$, endowed with the weak topology.
\end{theorem}

The proof can be found in Section~\ref{S:proofs}. It is an adaptation of
a standard coupling technique in propagation of chaos. We will couple our
process with interaction with an i.i.d.\ sequence of non-linear CMJ trees by
using the same source of noise (the variables $(\mathcal{P}_u,
\omega_u)_{u \in \mathcal{U}}$), but thinning the process with the
solution of \eqref{eq:main} instead of using the empirical age
distribution. 

A similar coupling has already been used for other types of
population models, for instance in \cite{chevallier2017,
forien2025stochastic, calvez2022dynamics, henry2023time}. Our
approach bears some differences that we want to highlight. (i) At the
conceptual level, we carry the coupling out for trees and study the limit
of the empirical measure of genealogical quantities (beyond the ancestral
lineages already considered in \cite{calvez2022dynamics, henry2023time}). (ii)
At the technical level, a CMJ process cannot be represented as the
solution to a Poisson-driven stochastic differential equation (conversely
to the models in the references above). In particular, we cannot rely on
standard stochastic analysis techniques to control the distance between
the coupled processes. This requires us to develop new arguments, and a
key step in our proof is to compare the differences in the coupling with
a CMJ process with immigration, which is a process that we introduce to
replace Gr\"onwall's inequality for Poisson stochastic differential
equations. 

Let us close this section by illustrating our main result with two
corollaries. First, we compute the limit of the age structure of the
population, viewed as a measure-valued stochastic process.

\begin{corollary} \label{cor:ageStructure}
    Under Assumption~\ref{ass:PDE}, if $(u_t)_{t \ge 0}$ is the solution
    to \eqref{eq:main},
    \[
        \lim_{N \to \infty} \big(\mu^N_t(\diff a);\, t \ge 0\big) 
        =
        (u_t(a)\diff a;\, t \ge 0),
    \]
    in probability for the topology of uniform convergence on compact
    sets.
\end{corollary}

The second corollary is concerned with a more genealogical quantity. For
an individual $u \in \mathcal{U}_i$ descended from the $i$-th ancestor of the
population, we can record a finite sequence of random variables giving
the successive birth times of its parents, namely $(\sigma_v)_{i \preceq
v \preceq u}$. We view it as a finite sequence of length $\abs{u}$
indexed backward-in-time, that is with first element $\sigma_u$ and last
element $\sigma_i$. 

The limit of the empirical measure of birth chains will be expressed 
in terms of the following Markov chain $(T_j)_{j \ge 1}$ on $\R$.
Conditional on $T_i = t > 0$, the chain jumps to a smaller time
distributed as 
\begin{equation} \label{eq:defRenewal}
    \E[\phi(T_{i+1}) \mid T_i = t ] = \frac{C(t,u_t)}{u_t(0)} 
    \int_0^\infty \phi(t-a) u_t(a) \tau(a) \diff a,
\end{equation}
where $(u_t(a);\, t,a \ge 0)$ is the solution to \eqref{eq:main}.
Setting $\phi \equiv 1$, the boundary condition of \eqref{eq:main}
implies that \eqref{eq:defRenewal} does define a probability distribution.
We stop the chain upon reaching $(-\infty, 0]$ and think of $(T_j)_j$
as being the finite collection of times until it hits the negative
half-line. The following result indicates that, conditional on observing
a birth at time $t > 0$, the birth chain leading to that
focal birth is asymptotically distributed as $(T_j)_j$, started from $T_1 =
t$.

\begin{corollary} \label{cor:infectionChain}
    If Assumption~\ref{ass:PDE} holds then, for any $T > 0$,
    \[
        \lim_{N \to \infty}
        \frac{1}{N} \sum_{i=1}^N \sum_{u \in \mathcal{U}_i} 
        \indic_{\{ \sigma^N_u \le T \}} \delta_{ (\sigma_v)_{i \preceq v \preceq u} }
        =
        \int_{-\infty}^T u_t(0) \mathscr{L}\big( (T_j)_j \mid T_1 = t \big) \diff t,
    \]
    in law as random measures on $\bigcup_{n \ge 1} \R^n$, endowed with
    the weak topology.
\end{corollary}

The limiting process has a simple intuitive description. In the limit,
there are $u_t(0)$ births at time $t$. Moreover, the expected number of
births at time $t$ generated by individual of age $a$ is $\tau(a) u_t(a)
C(t,u_t)$. Therefore, in expectation, the probability that the parent of
an individual born at time $t$ was born at time $s$ is proportional to
$\tau(t-s)u_t(t-s)$. The previous result shows that this simple
description becomes asymptotically exact as $N \to \infty$.

Finally, let us point at the fact that these are only two natural
functionals of the forest whose convergence follows from
Theorem~\ref{thm:main}. Other interesting examples can be found for
instance in \cite{Nerman1984, Coste2021}.

\subsection{Application to epidemiology}

In this section, we provide a brief epidemiological interpretation of our
results. Before going into this, we first mention that  our approach
provides an extension of \cite{Barbour2013, duchamps2023general}, who
consider a modification of a CMJ process with a fixed population size,
thereby introducing interactions among individuals. Their scaling limit
is derived by analysing the infection graph, a method that crucially
relies on the assumption of a fixed population size and on a specific
choice of the contact rate. In contrast, our approach—based on
propagation of chaos—is conceptually different and more direct, allowing
us to handle a broader class of population interactions.

In an epidemiological context, births correspond to infection events, and
the function $C$ represents the contact rate within the general
population. This rate may vary due to the implementation of control
measures or behavioral changes, in which case a fraction of potential
infections is prevented. We model this effect by pruning the
corresponding branches of the infection tree. Our formalism is fairly
general, allowing the contact rate to depend not only on the current
state of the epidemic but also on its past trajectory. We illustrate this
with two simple examples below.

\begin{example}
    Let $C$ be given by
    \[
        C(t, \mu) = 1 - \frac{\abs{\mu}}{K}
    \]
    where $K > 1$ and $\abs{\mu}$ is the total mass of $\mu$. When $k$
    individuals have been infected in the population, a potential
    infection is only effective with probability $1 - k/(NK)$. This
    models a closed population of fixed size $KN$ in which individuals
    become perfectly immune to the disease after infection. It is the
    same model as that considered in \cite{duchamps2023general}.
\end{example}

\begin{example}
    One can extend the previous example by letting 
    \[
        C(t, \mu) = \Big( 1 - \frac{\abs{\mu}}{K} \Big) \big(1 -
        \kappa \indic_{\{\abs{\mu} > \theta K\}} \big).
    \]
    The additional term models a control measure which is implemented
    when a fraction $\theta$ of the population has been infected, and
    which reduces the contact rate in the population to $1-\kappa$.
\end{example}

In this work, the CMJ process has been described solely through its
branching structure, as we only record the birth (or infection) times of
individuals. In particular, we do not model explicitly an individual's
death (or recovery). A natural extension consists in enriching the model
by associating to each individual not only a birth time but also a
stochastic process describing its evolution over time.  

Formally, a CMJ process can then be defined from a pair $(\mathcal{P},
X)$, where $\mathcal{P}$ is a point process and $X \equiv (X(a))_{a \ge 0}$
is a stochastic process with values in a discrete state space $S$
modeling observable individual quantities such as clinical status
(symptom onset, hospitalisation, recovery); see \cite{foutel2022individual, 
duchamps2023general} for concrete examples in epidemiology. The process
$X$ is an example of a \emph{random characteristic}
\cite{jagers1975branching, Nerman1981}.

While this additional structure is highly relevant from a modeling
perspective, it primarily complicates notation without raising new
mathematical difficulties. For this reason, we did not include
characteristics in the main text. However, our results extend directly to
the framework of \cite{foutel2022individual} with minor additional
bookkeeping. For instance, in the spirit of \cite{foutel2022individual,
duchamps2023general}, Corollary~\ref{cor:ageStructure} can be extended to
describe the limiting behavior of epidemic observables:
\[
    \forall x \in S, \quad 
    \frac{1}{N}\sum_{i=1}^N \sum_{u \in {\cal U}_i} \mathbf{1}_{\{\sigma^N_u \le t\}}  
    \mathbf{1}_{\{X_u(t-\sigma^N_u) = x\}}
    \;\Longrightarrow\;
    \int_0^\infty \mathbb{P}(X(a) = x) u_t(a) \diff a.
\]

Finally, in epidemiology, an important quantity of interest is the
reconstruction, backward in time, of the chain of infections that led to
the infection of a focal individual. Informally,
Corollary~\ref{cor:infectionChain} states that, conditional on an
infection occurring at time $t$, this infection chain follows a
time-inhomogeneous Markov chain.

\subsection{Outline}

The rest of the manuscript is organised as follows. Section~\ref{S:PDE}
recalls the definition of a solution to \eqref{eq:main} from
\cite{webb1985theory} and proves that this equation is well-posed. As indicated
above, our main result will follow from a coupling argument.
Section~\ref{S:CMJImmigration} introduces a new process, the CMJ process
with immigration, which is used to dominate the error made in the
coupling. The bulk of the technical work is carried out in
Section~\ref{S:coupling}, where the coupling is defined and is compared
to an appropriate CMJ process with immigration. Finally, the proofs are
completed in Section~\ref{S:proofs}, where we also provide the necessary
continuity results to derive our two corollaries.

\section{Well-posedness of the limit}
\label{S:PDE}

We define an appropriate notion of weak solution to the PDE, which is based on
the monograph \cite{webb1985theory}. Let $g \in \mathrm{L}^1(\R_+)$. We
say that a function $(u_t(a);\, t,a \ge 0)$ is a solution to 
\begin{align*}
    (\partial_t + \partial_a) u_t(a) &= 0 \\
    u_t(0) &= C\big(t, u_t\big) \int_0^\infty u_t(a) \tau(a) \diff a\\
    u_0(a) &= g(a)
\end{align*}
if, for all $t \ge 0$, $u_t \in \mathrm{L}^1(\R_+)$ is such that for a.e.\ $a
\in [0,t)$
\[
    u_t(a) = C(t-a, u_{t-a}) \int_0^\infty
    u_{t-a}(s) \tau(s) \diff s,
\]
and $u_t(a) = g(a-t)$ for a.e.\ $a \ge t$. This definition is motivated
by a formal application of the method of characteristics which suggests
that a (strong) solution to \eqref{eq:main} should satisfy that 
\[
    u_t(a) = u_{t-a}(0) = C(t,u_{t-a}) \int_0^\infty u_{t-a}(s) \tau(s) \diff s,
\]
for $a \le t$, and $u_t(a) = u_0(a-t) = g(a-t)$ for $a \ge t$. We refer
to \cite{webb1985theory} for more details on this type of age-structured
equations.

\begin{proposition}
    For any $g \in \mathrm{L}^1(\R_+)$, under Assumption~\ref{ass:PDE}, there
    exists a unique solution to \eqref{eq:main} with initial condition
    $g$.
\end{proposition}

\begin{proof}
    We cannot rely directly on the general results of \cite{webb1985theory} 
    because of our time inhomogeneity but the proof -- which relies
    on standard arguments -- can be readily adapted. Let us provide some
    details, since our setup is less general than that of
    \cite{webb1985theory} and allows for several simplifications. We
    start by showing local existence. Namely we prove that, for any $r >
    0$, there is $T \equiv T(r)$ such that \eqref{eq:main} has a unique
    solution on $[0,T]$ when started from any initial condition $g$ with
    $\norm{g}_1 \le r$.

    For $t \ge 0$ and $u \equiv (u_t)_{t \le T}$, let $K_t[u] \in
    \mathrm{L}^1(\R_+)$ be defined as
    \[
        K_t[u](a) = C(t-a, u_{t-a}) \int_0^\infty u_{t-a}(s) \tau(s) \diff s,
    \]
    for a.e.\ $a \in [0,t)$ and $K_t[u](a) = u_0(t-a)$ for a.e.\ $a \ge
    t$. We will prove the result by showing that $K \colon (u_t)_{t \le
    T} \mapsto (K_t[u])_{t \le T}$ is a contraction on 
    \[
        M \coloneqq \big\{ u \colon [0,T] \to \mathrm{L}_1(\R_+) \text{
            continuous, } \sup_{t \le T} \norm{u_t}_1
        \le 2r,\: u_0 = g \big\},
    \]
    for the uniform norm and some well-chosen $T > 0$. Let us write
    $\angle{f,g}$ for $\int f(a)g(a) \diff a$. For $(u_t)_{t\le T} \in
    M$ and $t \le T$,
    \begin{align*}
        \norm{K_t[u]}_1 
        &= \int_0^t C(t-a, u_{t-a}) \angle{u_{t-a}, \tau} \diff a + \norm{g}_1 \\
        &\le 2r \norm{\tau}_\infty T + r.
    \end{align*}
    Moreover, for $t \le t' \le T$,
    \begin{multline*}
        \norm{K_t[u] -K_{t'}[u]}_1
        \le 
        \int_{t'}^\infty \abs{g(a-t) - g(a-t')} \diff a + \int_0^{\abs{t-t'}} g(a) \diff a \\
        + 2 r \abs{t-t'} \norm{\tau}_\infty
        + \int_0^t \abs[\big]{C(t-a, u_{t-a}) \angle{u_{t-a}, \tau} - C(t'-a,
            u_{t'-a}) \angle{u_{t'-a}, \tau}} \diff a.
    \end{multline*}
    All four terms vanish as $t - t' \to 0$ (by continuity of the shift
    in $\mathrm{L}^1(\R_+)$ for the first term, and dominated convergence for
    the last one), showing that $K[u] \colon [0,T] \to \mathrm{L}^1(\R_+)$ is
    continuous. Thus, $K$ maps $M$ to itself for $T$ small enough.

    We now show that $K$ is a contraction. Let $u, u' \in M$,
    \begin{align*}
        \norm{K_t[u] - K_t[u']}_1
        &\le \begin{multlined}[t]
        \int_0^t \Big( \abs[\big]{C(t-a, u_{t-a}) -C(t-a, u'_{t-a})}
        \int_0^\infty u_{t-a}(s) \tau(s)\diff s \\
        + C(t-a, u'_{t-a}) \int_0^\infty \abs{u_{t-a}(s) - u'_{t-a}(s)} \tau(s)\diff s
        \Big) \diff a
        \end{multlined} \\
        &\le T \big(L 2r \norm{\tau}_\infty + 1 \big) 
            \sup_{t \le T} \norm{u_t - u'_t}_1.
    \end{align*}
    By choosing $T$ small enough, $K$ is a strict contraction and has a
    unique fixed point in $M$, which is a solution to \eqref{eq:main} on
    $[0,T]$.

    Global existence will follow from local existence and an \emph{a
    priori} bound. If $(u_t)_{t \ge 0}$ is a solution to \eqref{eq:main}, 
    \begin{align*}
        \forall t \ge 0,\quad 
        \norm{u_t}_1 
        &= \int_0^t C(t-a,u_{t-a}) \int_0^\infty u_{t-a}(s) \tau(s) \diff s \,\diff a 
        + \norm{u_0}_1 \\
        & \le \norm{\tau}_\infty \int_0^t \norm{u_{t-a}}_1 \diff a +
        \norm{u_0}_1.
    \end{align*}
    By Gr\"onwall's inequality, 
    \[
        \forall t \ge 0,\quad \norm{u_t}_1 \le \norm{u_0}_1 e^{\norm{\tau}_\infty t}.
    \]
    Fix $t > 0$. We can construct a solution to \eqref{eq:main} on
    $[0,t]$ by covering this interval with finitely many subintervals of
    length $T(r)$ for $r = \norm{g}_1 e^{\norm{\tau}_\infty t}$, and
    constructing the solution on each subinterval using the first part of
    the proof. This solution is unique by another application of
    Gr\"onwall's inequality, and it can be extended to a unique solution
    on the whole half line.
\end{proof}

\begin{proof}[Proof of Proposition~\ref{prop:nonLinearExistence}]
    Suppose that \eqref{eq:main} has a solution $(u_t)_{t \ge 0}$ and let
    $\widetilde{\mu}_t$ be the empirical distribution of ages of a CMJ
    tree started from one ancestor and thinned by the time-varying
    function $(C(t, u_t))_{t \ge 0}$. Then \cite[Corollary~7]{foutel2022individual} 
    shows that 
    \[
        \E[\widetilde{\mu}_t(\diff a)] = v_t(a) \diff a,
    \]
    where $(v_t(a); t,a \ge 0)$ is the unique solution to
    \begin{align}
    \begin{split} \label{eq:McKendrickVonFoster}
        (\partial_t + \partial_a) v_t(a) &= 0 \\
        v_t(0) &= C(t, u_t) \int_0^\infty v_t(a) \tau(a) \diff a\\
        v_0(a) &= g(a).
    \end{split}
    \end{align}
    Since $u_t(a) \diff a$ is also solution to this linear equation one
    must have that $v_t(a) \diff a = u_t(a) \diff a$. Therefore the
    corresponding CMJ tree is a non-linear tree.

    Conversely, if there exists a non-linear CMJ tree with empirical
    age distribution $\mu'_t$, then $u_t(a) \diff a = \E[\mu'_t(\diff
    a)]$ solves \eqref{eq:McKendrickVonFoster} and thus provides a
    solution to \eqref{eq:main}. 
\end{proof}

\section{Crump--Mode--Jagers process with immigration}
\label{S:CMJImmigration}

\begin{figure}[t]
    \centering
    \includegraphics[width=\textwidth]{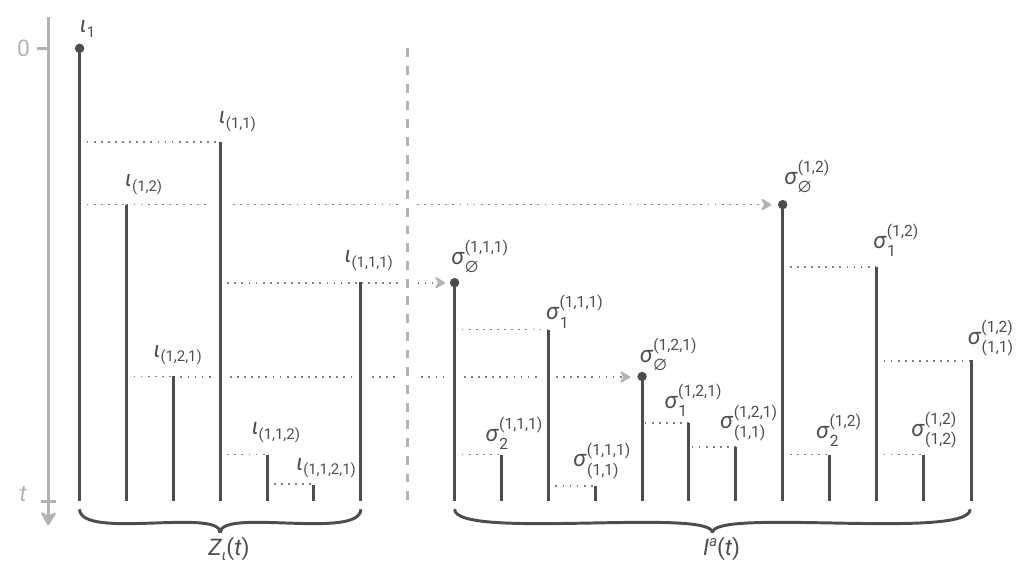}
    \caption{
        Illustration of the CMJ process with immigration.
        The left panel corresponds to the CMJ tree of potential
        immigration times, which is started from a unique ancestor. At
        each potential immigration time, a new particle might be added to
        the right panel (represented by an arrow) and starts to reproduce
        as a CMJ process. The process $(I^\eta(t))_{t \ge 0}$ counts the
        number of particle on the right.
    }
    \label{fig:cmj_immigration}
\end{figure}

This section introduces a Crump--Mode--Jagers process in which new
individuals are added to the population at random times. For our purpose,
these immigration times are obtained by thinning the set of birth times
of a different, independent CMJ process. The definition of this process is
illustrated in Figure~\ref{fig:cmj_immigration}. As highlighted in the introduction, this process will be used in later sections to dominate the error made in the
coupling between the non-linear CMJ process and the CMJ process with immigration. 

For $i \ge 1$, let $( \iota_u )_{u \in \mathcal{U}_i} \in \mathscr{T}$ be
independent CMJ trees (without thinning) as described at the beginning of
Section~\ref{SS:constructionCMJ}. In particular, the initial ages
$(\iota_i)_{i \ge 1}$ of the individuals in the first generation are
distributed as $g(a)\diff a$. These birth times correspond to the set of
potential immigration times.
Let $(\omega_u)_{u \in \mathcal{U}}$ be i.i.d.\ uniform random variables
on $(0,1)$. For each $u \in \mathcal{U}$, let $(\sigma^u_v)_{v \in
\mathcal{U}}$ be an independent CMJ tree, started from one ancestor with
birth time $\sigma^u_{\varnothing} = 0$. (Note that we use a slightly
different labeling compared to Section~\ref{SS:constructionCMJ} and let
$\varnothing$ be the ancestor, which is the standard convention.)
Further, let $Z^u(t)$ be the population size in the CMJ tree
$(\sigma^u_v)_{v \in \mathcal{U}}$, that is,
\[
    \forall t \ge 0,\quad Z^u(t) = \sum_{v \in \mathcal{U}} \indic_{\{ \sigma^u_v \le t \}}.
\]
We insist on the fact that we attach an independent CMJ process
$(Z^u(t))_{t \ge 0}$ to \emph{each} individual $u$ in $\mathcal{U}$.

Fix $\eta > 0$ and $N \ge 1$. Recall that $L$ is the Lipschitz constant of Assumption \ref{ass:PDE}. We construct a process
$(I^\eta(t))_{t \ge 0}$, called a \emph{CMJ process with immigration}, as
\[
    \forall t \ge 0,\quad
    I^\eta(t) = \sum_{i=1}^N \sum_{u \in \mathcal{U}_i} 
    \indic_{\{ 0 < \iota_u \le t\}}
    \indic_{\{ \omega_u \le \eta + \frac{L}{N}I^\eta(\iota_u^-)\}}
    Z^u(t-\iota_u).
\]
The process $(I^\eta(t))_{t \ge 0}$ can be obtained from a particle system
with two types. The system starts with $N$ particles of type $1$ with
initial ages distributed as $g(a)\diff a$. Particles of type $1$
reproduce independently according to the CMJ dynamics described in
Section~\ref{SS:constructionCMJ}, and correspond to the left of
Figure~\ref{fig:cmj_immigration}. At each birth of a particle of type $1$
(say at time $t$), a new particle of type $2$ is added (immigrated) with
probability $(\eta + LI^\eta(t^-)/N) \wedge 1$. Once immigrated, particles of
type $2$ follow the CMJ dynamics of Section~\ref{SS:constructionCMJ}
independently. The process $(I^\eta(t))_{t \ge 0}$ counts the number of type
$2$ particles.

We emphasise that we use a slight abuse of notation when referring to
this process as a CMJ process with immigration. Conversely to the usual
setting of branching processes with immigration, potential immigration
times in the right population (in Figure~\ref{fig:cmj_immigration}) are
not Poisson, but subordinated to the birth times of the left population
(in Figure~\ref{fig:cmj_immigration}), which are the birth times of a CMJ
process. Moreover, an immigration is only `successful' with a
probability that depends on the number of particles in the right
population, introducing a dependence between the state of the system and
the immigration times.

Note that this process has a self-exciting behaviour: the more particles
in the population, the higher the probability of immigrating a new
particle. The following result will provide enough control on the number
of particles for our purpose.

\begin{lemma} \label{lem:convImmigration}
    Fix some $t > 0$ and $L > 0$. For any $\epsilon > 0$,
    \[
        \adjustlimits \lim_{\eta \to 0} \limsup_{N \to \infty} 
        \P(I^\eta(t) \ge \epsilon N)
        = 0.
    \]
\end{lemma}

\begin{proof}
    Let us find a first upper bound for $I^\eta(t)$. Define $S^\eta(s)$ as 
    \begin{equation} \label{eq:boundBranching}
        S^\eta(s) = \sum_{i=1}^N 
        \sum_{u \in \mathcal{U}_i} \indic_{\{ 0 < \iota_u \le s \}} 
        \indic_{\{ \omega_u \le \eta + \tfrac{L}{N} S^\eta(\iota_u^-) \}}
        Z^u(t).
    \end{equation}
    Recall that $Z^u(t)$ records the number of births up to time $t$ in
    the CMJ originated from $u$. Since $Z^u(t) \ge Z^u(s-\iota_u)$ as
    long as $s \le t$, an induction on the birth times shows that $S^\eta(s)
    \ge I^\eta(s)$ as long as $s \le t$.

    Let $\widetilde{\iota}_1 \le \widetilde{\iota}_2 \le \dots$ be
    obtained by re-ordering the \emph{positive} birth times of the
    immigration process, namely re-ordering the set
    \[
        \bigcup_{i=1}^N \big\{ \iota_u : u \in \mathcal{U}_i, u \ne i \big\}.
    \]
    Set $\iota_0 = 0$. Then, if $S^\eta_n = S^\eta(\widetilde{\iota}_n)$, it
    follows from \eqref{eq:boundBranching} that $(S^\eta_n;\, n \ge 1)$ is a
    Markov chain such that, conditional on $S^\eta_n$, 
    \[
        S^\eta_{n+1} =  
        \begin{cases}
            S^\eta_n + Z_n &\text{with probability
            $(\eta+\tfrac{L}{N}S^\eta_n) \wedge 1$} \\
            S^\eta_n &\text{otherwise,}
        \end{cases}
    \]
    where $(Z_n)_{n \ge 1}$ are i.i.d.\ with $Z_n \sim Z(t)$. Therefore,
    \[
        \E[S^\eta_{n+1}] = \E[S^\eta_n] + \E[Z(t)]
        \E[(\eta+\tfrac{L}{N}S^\eta_n) \wedge 1].
    \]
    If we let $(x_n)_{n \ge 1}$ be the solution to 
    \[
        x_{n+1} = x_n \big( 1 + \tfrac{L}{N} \E[Z(t)] \big) + \eta \E[Z(t)], 
        \qquad x_0 = 0,
    \]
    a simple induction shows that $\E[S^\eta_n] \le x_n$.
    Solving this recursion explicitly yields that
    \[
        \E[S^\eta_n] \le x_n = 
        \frac{\eta N}{L} \Big( \big(1+\tfrac{L}{N}\E[Z(t)]\big)^n - 1\Big).
    \]

    To conclude the proof, let $(Z_\iota(t))_{t \ge 0}$ be the process
    that counts the size of the immigration process, namely
    \[
        Z_\iota(t) = \sum_{i=1}^N \sum_{u \in \mathcal{U}_i} 
        \indic_{\{ \iota_u \le t\}}.
    \]
    Note that $S^\eta(t) = S^\eta_{Z_\iota(t)-N} \le S^\eta_{Z_\iota(t)}$. 
    (The offset by $N$ is due to the fact the birth times of the $N$
    ancestors are not potential immigration times.) By the law of large
    numbers, $Z_\iota(t) / N \to \E[Z_1(t)]$ almost surely, where
    $(Z_1(t))_{t \ge 0}$ is a CMJ process started from one ancestor with
    age distributed as $g(a)\diff a$. Therefore, for $K > \E[Z_1(t)]$,
    Markov's inequality entails that 
    \begin{align*}
        \P( S^\eta(t) \ge \epsilon N ) 
        &\le \P( Z_\iota(t) \ge K N) + \P( S^\eta_{KN} \ge \epsilon N ) \\
        &\le \P( Z_\iota(t) \ge K N) + \frac{\eta}{L\epsilon} 
            \Big[ 
                \big( 1 + \tfrac{L}{N} \E[Z(t)] \big)^{NK} - 1
            \Big]  \\
        &\le \P( Z_\iota(t) \ge K N) + 
            \frac{\eta}{L\epsilon} \big( e^{KL\E[Z(t)]} - 1 \big).
    \end{align*}
    The result follows by letting $N \to \infty$ first, and then $\eta \to 0$.
\end{proof}

\section{Construction of the coupling}
\label{S:coupling}

\subsection{Algorithmic construction of CMJ processes}
\label{SS:algoCMJ1}

Our coupling with a CMJ process with immigration necessarily requires
cumbersome bookkeeping. It relies on an `algorithmic' construction of
CMJ processes with thinning, in which potential birth times are inspected
in increasing order and the population is updated depending on whether
the birth is effective or is thinned. The coupling procedure is rather involved, and as a warm up, we first illustrate this
construction in the case of the CMJ process with interaction, and carry
out the full coupling in the next section.

Let $(\mathcal{P}_n)_{n \ge 1}$ and $(\omega_n)_{n \ge 1}$ be independent
i.i.d.\ sequences, distributed as $\mathcal{P}$ and uniformly on $[0,1]$,
respectively. We construct iteratively two sequences of measures
on $\R$: $L^N_n$ and $R^N_n$. Informally, the $n$-th step of the
algorithm inspects a potential birth in the CMJ process, at time
$\sigma^N_n$. The atoms of the measure $R^N_n$ (for `right') record the
future birth times to be inspected. This is the `coming generation' in
CMJ terminology, which corresponds to the individuals not born at time
$\sigma^N_n$, but whose parents are already born at that time.
Conversely, $L^N_n$ (for `left') records the birth times of individuals
which are born before time $\sigma^N_n$, and have not been thinned. It
encodes the state of the population at time $\sigma^N_n$. Let us define
the variables rigorously. We drop the superscript $N$ to ease the
notation.

Let $(A^0_i, \mathcal{P}^0_i)_{i \ge 1}$ be i.i.d.\ random pairs
distributed as in \eqref{eq:initialCondition}. We initialise the
variables to 
\[
    L_0 = \sum_{i=1}^N \delta_{-A^0_i},\qquad
    R_0 = \sum_{i=1}^N \sum_{a \in \mathcal{P}^0_i} \delta_a,\qquad
    \sigma_{-1} = 0.
\]
Suppose that $(L_n, R_n)$ have been constructed, and let $\sigma_n$ be
the smallest atom of $R_n$. Define the empirical age distribution 
at step $n$ (that is at time $\sigma_n^-$) as
\[
    \forall n \ge 0,\quad \mu_{\sigma^-_n} = \frac{1}{N} \sum_{s \in L_n}
    \delta_{\sigma_{n}-s}.
\]
On the event $\{ \omega_n \le C(\sigma_n, \mu_{\sigma_n^-})\}$ the birth
is effective in which case 1) an atom is added to $L_n$ that record the
new birth at time $\sigma_n$ and 2) a copy of $\mathcal{P}$ shifted by
$\sigma_n$ is added to $R_n$, that records the future times at which
this individual gives potential birth.

More formally, we define
\[
    L_{n+1} = L_n + 
    \indic{\{ \omega_n \le C(\sigma_n, \mu_{\sigma^-_n}) \}} \delta_{\sigma_n}
    ,\quad
    R_{n+1} = R_n - \delta_{\sigma_n}
    + \indic{\{ \omega_n \le C(\sigma_n, \mu_{\sigma^-_n}) \}} 
    \sum_{a \in \mathcal{P}_n} \delta_{\sigma_n+a}.
\]
If $R_n$ has no atoms, we set $R_{n+1} = 0$, $\sigma_n = \infty$, and
$L_{n+1}=L_n$.

We can define a time-dependent empirical measure by setting
\begin{equation}
\label{eq:mu_n:mu_t}
    \mu_t = \mu_{\sigma_n^-}, \qquad \text{for $t \in [\sigma_{n-1}, \sigma_n)$.} 
\end{equation}
It should be clear that $(\mu_t)_{t \ge 0}$ is distributed as the
empirical distribution of ages in the CMJ process with interaction
constructed in Section~\ref{SS:constructionCMJ}:
a potential birth at time $t$ is added to $L_n$ with probability $C(t,
\mu_{t^-})$ and discarded otherwise. On the event that the birth is
effective, an independent progeny is added to $R_n$, which is distributed
as the atoms of $\mathcal{P}$ shifted by the birth time $\sigma_n$ of the
mother.

\subsection{Coupling with immigration}
\label{SS:fullCoupling}

The aim of the section is to prove the following result. It provides a
coupling of the CMJ tree with interaction to an i.i.d.\ collection of
non-linear CMJ trees, with the difference between the two processes in
the coupling being controlled by a CMJ process with immigration as described in the previous section.

\begin{proposition} \label{prop:coupling}
    Suppose that Assumption~\ref{ass:PDE} holds.
    For any $a > 0$ and $N \ge 1$, there exists a coupling between 
    \begin{itemize}
        \item a collection of CMJ trees with interaction 
            $(\mathcal{T}^N_i)_{i \le N}$ started from $N$ individuals;
        \item i.i.d.\ copies of non-linear CMJ trees $(\mathcal{T}'_i)_{i
            \le N}$; and
        \item a CMJ process with immigration $(I^a(t))_{t \ge 0}$.
    \end{itemize}
    Let the empirical age distribution of the non-linear trees be 
    \[
        \forall t \ge 0,\quad \overline{\mu}_t^N(\diff a) 
        = \frac{1}{N} \sum_{i=1}^N \sum_{u \in \mathcal{U}_i} 
        \indic_{\{\sigma'_u \le t\}}\delta_{t-\sigma'_u}(\diff a).
    \]
    The coupling can be chosen so that, on the event $\{ L
    d_{\Pr}(\overline{\mu}^N_T, u_T) < \eta\}$,
    \begin{equation} \label{eq:goodCoupling}
        \sum_{i=1}^N \sum_{u \in \mathcal{U}_i} 
        \indic_{\{\sigma^N_u \le t, \sigma'_u \le t\}} 
        \indic_{\{ \sigma^N_u \ne \sigma'_u\}} 
        \le I^\eta(t).
    \end{equation}
\end{proposition}

Let us provide a verbal description of the coupling before providing a
formal proof.

\paragraph{Coupling the two thinning procedures.}
In order to carry the coupling out, a potential birth at time $t$ 
needs to be kept with probability $C(t, \mu_{t^-})$ in the process with
interaction, and with probability $C(t, u_t)$ in the non-linear process.
It will considerably ease the discussion to introduce the notation
\[
    P_1(t) = C(t, \mu_{t^-}),\qquad P_2(t) = C(t, u_t)
\]
for these two probabilities. 

We couple the two procedures by adding labels to the particles that
record for which procedure they have been thinned. There are three
labels: $(1,1)$, $(1,0)$, and $(0,1)$.
A zero in the first (resp.\ the second) coordinate is interpreted as
saying that the particle has been thinned with probability $P_1(t)$ 
(resp.\ $P_2(t)$). Hence, particles $(1,1)$ are present for both
procedures, whereas particles $(1,0)$ and $(0,1)$ have been thinned for
one procedure but not for the other. We need to control the number of
particles of type $(1,0)$ and $(0,1)$.

Initially, all particles are present in both processes and have label
$(1,1)$. Consider a particle born of a $(1,1)$ parent at time $t$. The
label $\ell$ of this new particle is assigned in such a way that the
probability that it has a $1$ in the first (resp.\ the second) column is
$P_1(t)$ (resp.\ $P_2(t)$). We do so in a way that maximises the
probability that the particle is well-coupled (that is of type $(1,1)$)
by setting
\begin{equation} \label{eq:lawCoupling}
    \P\big( \ell = (i,j) \big) = 
    \begin{dcases}
        P_1(t) \wedge P_2(t) &\text{if $(i,j)=(1,1)$,}\\
        (P_1(t) - P_2(t)) \vee 0 &\text{if $(i,j)=(1,0)$,}\\
        (P_2(t) - P_1(t)) \vee 0 &\text{if $(i,j)=(0,1)$.}\\
    \end{dcases}
\end{equation}
Indeed, $\P(\ell \in \{(1,1), (1,0)\}) = P_1(t)$ and $\P(\ell \in
\{(1,1), (0,1)\}) = P_2(t)$.

\begin{figure}[t]
    \centering
    \includegraphics[width=\textwidth]{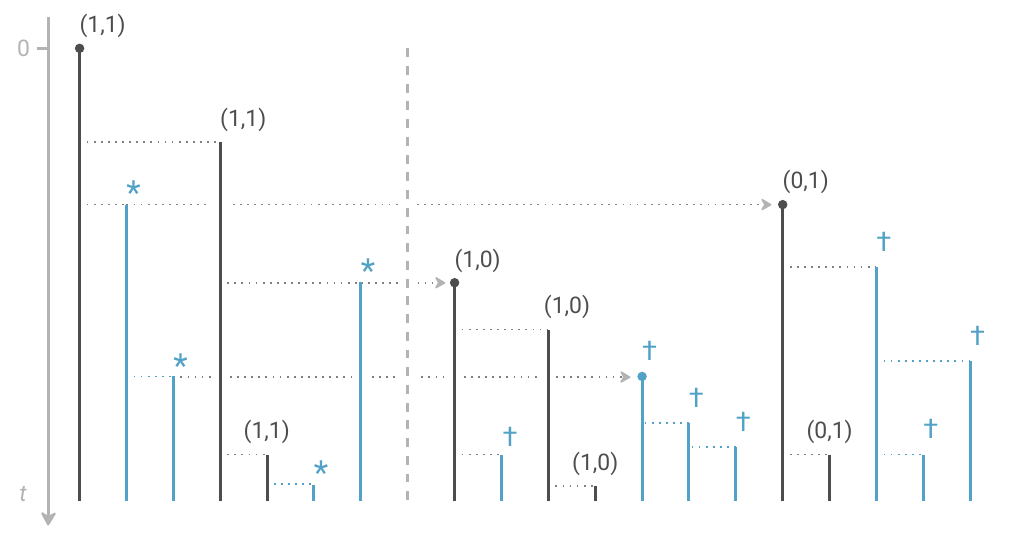}
    \caption{Illustration of the coupling, starting from a single
        particle. Non-fictitious individuals in the coupling (with types
        $(1,1)$, $(1,0)$, and $(0,1)$) are represented by dark lines,
        whereas fictitious individuals (with types $\dag$ and $*$) by
        blue lines. The birth times of individuals on the left (with
        types $(1,1)$ and $\dag$) correspond to the potential immigration
        times in the CMJ process with immigration. The tree of $(1,1)$
        particles is thinned, but it is complemented with $*$ particles
        so that it is distributed as CMJ tree. Whenever a particle is
        immigrated and added to the right of the picture, it produces a
        tree which we complement with $\dag$ particles so that it is
        also distributed as a CMJ tree.}
    \label{fig:coupling}
\end{figure}

\paragraph{Comparison with a CMJ process with immigration.} Although the
previous construction provides a natural coupling 
between the process with
interaction and independent non-linear processes, controlling the number
of particles of type $(1,0)$ and $(0,1)$ is not an easy task. We do so by
augmenting the process with two new types of particles ($*$ and $\dag$)
to construct a CMJ process with immigration dominating the discrepancy
between the two processes. 
This is illustrated on Fig. \ref{fig:coupling} and can be described as follows.

Our objective is to control the size of the population of particles $(1,0)$
and $(0,1)$. A $(1,0)$ (resp.\ $(0,1)$) particle reproduces as a CMJ
process thinned with probability $P_1(t)$ (resp.\ $P_2(t)$). The tree
grown out of such a particle can be dominated by an un-thinned CMJ
process by adding a particle (of type $\dag$) each time the offspring of
a $(1,0)$ (or $(0,1)$) particle is thinned. The particles $\dag$
correspond to fictitious individuals and we let a $\dag$ particle reproduce
as an un-thinned CMJ process of type $\dag$. 

A new particle of type $(1,0)$ or $(0,1)$ is produced with probability
$\abs{P_1(t) - P_2(t)}$ whenever a $(1,1)$ individual gives birth at time
$t$. We envision the addition of such a particle as an `immigration' (right part of Fig. \ref{fig:coupling}).
At first, the set of all potential immigration times is the set of birth
times from a $(1,1)$ parent. This tree of $(1,1)$ individuals can also be
dominated by an un-thinned CMJ process by adding to the population a
particle of type $*$ whenever the offspring of a $(1,1)$ individual is
not of type $(1,1)$, and letting $*$ particles reproduce as an un-thinned
CMJ process (see left part of Figure~\ref{fig:coupling}). Finally, the
probability of immigrating a new $(1,0)$ or $(0,1)$ particle at the birth
time $t$ of a $(1,1)$ particle is $\abs{P_1(t)-P_2(t)}$. Using the
Lipschitz condition of Assumption~\ref{ass:PDE}, we will show that 
\[
    \abs{P_1(t) - P_2(t)} \le \frac{L}{N} I^\eta(t) + \eta
\]
for $t \le T$ with high probability and where $I^\eta(t)$ counts the
number of births in the right part of Figure~\ref{fig:coupling} by time
$t$. The coupling will be completed by immigrating a $\dag$ particle with
probability $(\frac{L}{N} I^\eta(t) + \eta) \wedge 1 - \abs{P_1(t)-
P_2(t)}$ at each birth of a $(1,1)$ parent.

We now move on to the proof, which amounts to formalising the coupling
described above.

\begin{proof}[Proof of Proposition~\ref{prop:coupling}]
Consider an i.i.d.\ sequence $(\omega_n)_{n \ge 1}$ of uniform r.v.s on
$(0,1)$, and three additional i.i.d.\ sequences $(\mathcal{P}_n)_{n \ge
1}$, $(\mathcal{P}^*_n)_{n \ge 1}$, and $(\mathcal{P}^\dag_n)_{n \ge 1}$
distributed as $\mathcal{P}$. As in Section~\ref{SS:algoCMJ1}, we
construct two sequences of measures $(L_n)_{n \ge 0}$ and $(R_n)_{n \ge
0}$. (All random variables depend on $N$, but we make this dependence
implicit.) However, conversely to Section~\ref{SS:algoCMJ1}, each $L_n$ and
$R_n$ will be a measure on $\R_+ \times E$, where $E$ is the set of labels
\[
    E = \{ (1,1), (1,0), (0,1), \dag, * \}.
\]
The measures $L_n$ records the state of the population at time
$\sigma_n^-$: each of its atoms $(t, \ell) \in L_n$ encodes an individual
with label $\ell$ born at time $t$. Similarly, $R_n$ records the `coming
generation': an atom $(t, \ell) \in R_n$ corresponds to a potential
future birth at time $t$, from a \emph{parent} with label $\ell$.

It will be convenient to introduce the notation
\[
    \mu_{\sigma_n^-}(\diff s) 
    = \frac{1}{N} L_n\big(\diff s \times \{(1,1), (1,0)\}\big),
    \quad
    \overline{\mu}_{\sigma^-_n}(\diff s) 
    = \frac{1}{N} L_n\big(\diff s \times \{(1,1), (0,1)\}\big),
\]
for the empirical measures of birth times of the CMJ trees with
interaction and the non-linear trees, respectively. We will also need to
define
\[
    I^\eta({\sigma^-_n}) = 
    L_n\big( \R_+ \times \big\{ (1,0), (0,1), \dag \big\} \big).
\]
for some fixed $\eta > 0$. As in Section~\ref{SS:algoCMJ1} (see
\eqref{eq:mu_n:mu_t}), we can turn these sequences into real-time
processes by defining 
\[
    \forall t \in [\sigma_{n-1}, \sigma^-_n), \quad 
    I^\eta(t) = I^\eta(\sigma^-_n),\quad
    \mu_t = \mu_{\sigma^-_n},\quad
    \overline{\mu}_t = \mu_{\sigma^-_n},
\]
so that $I^\eta(t)$ captures the number of births on the right of 
Figure~\ref{fig:coupling}.

\medskip
\noindent
\textit{Initialisation.} Let $(A^0_i, \mathcal{P}^0_i)_{i \ge 1}$ be
i.i.d.\ as in \eqref{eq:initialCondition} and define  
\[
    L_0 = \sum_{i=1}^N \delta_{(-A^0_i, (1,1))}, \qquad
    R_0 = \sum_{i=1}^N \sum_{a \in \mathcal{P}^0_i} \delta_{(a,(1,1))},
    \qquad \sigma_{-1} = 0.
\]
Since the two processes start from the same initial condition, all
initial individuals are well-coupled with label $(1,1)$.
As above, let $(\sigma_n, \ell_n)$ be the atom of $R_n$ with the smallest
time coordinate. Recall that $\ell_n$ is the label of the parent giving
birth at time $\sigma_n$. We need to consider several cases, depending on
the value of $\ell_n$.

\medskip
\noindent \textit{Type $(1,0)$, $(0,1)$, and $\dag$ individuals.}
Suppose that $\ell_n \in \{(1,0), (0,1), \dag \}$.
Define $\ell'_n$ as 
\begin{equation} \label{eq:couplingDag}
    C_n =
    \begin{dcases}
        C(\sigma_n, \mu_{\sigma^-_n}) &\text{ if $\ell_n = (1,0)$,} \\
        C(\sigma_n, u_{\sigma_n}) &\text{ if $\ell_n  = (0,1)$,} \\
        1 &\text{ if $\ell_n = \dag$},
    \end{dcases}
    \qquad
    \ell'_n = 
    \begin{dcases}
        \ell_n &\text{if $\omega_n \le C_n$,} \\
        \dag &\text{else.}
    \end{dcases}
\end{equation}
The event $\{\ell'_n = \dag\}$ corresponds to thinning the offspring. We let 
\begin{equation*}
    L_{n+1} = L_n + \delta_{(\sigma_n, \ell'_n )},\qquad
    R_{n+1} = R_n - \delta_{(\sigma_n, \ell_n)} + 
    \sum_{a \in \mathcal{P}_n} \delta_{(\sigma_n+a, \ell'_n)}.
\end{equation*}

\medskip
\noindent \textit{Type $*$ individuals.}
Second, if $\ell_n = *$, let $\rho_n$ be the indicator of the event
$\{ \omega_n \le L I^\eta(\sigma^-_n)/N + \eta \}$. We define 
\begin{align}
\begin{split}\label{eq:couplingStar}
    L_{n+1} &= L_n + \delta_{(\sigma_n, *)} + \rho_n \delta_{(\sigma_n,
    \dag)}, \\
    R_{n+1} &= R_n - \delta_{(\sigma_n, *)} 
    + \sum_{a \in \mathcal{P}^*_n} \delta_{(\sigma_n+a,*)}
    + \rho_n \sum_{a \in \mathcal{P}_n} \delta_{(\sigma_n+a,\dag)}.
\end{split}
\end{align}
In words, $*$ individuals always produce one CMJ generation of new $*$
individuals and immigrate a new $\dag$ particle if $\rho_n = 1$.

\medskip
\noindent \textit{Type $(1,1)$ individuals.}
Finally, if $\ell_n = (1,1)$, define
\begin{equation} \label{eq:coupling11}
    \ell'_n = 
    \begin{dcases}
        (1, 0) &\text{if $\:\omega_n \le C(\sigma_n, \mu_{\sigma^-_n}) -
            C(\sigma_n, u_{\sigma_n})$,} \\
        (0, 1) &\text{if $\:\omega_n \le C(\sigma_n, u_{\sigma_n}) -
            C(\sigma_n, \mu_{\sigma^-_n})$,} \\
        (1, 1) &\text{if $\:\abs{C(\sigma_n, \mu_{\sigma^-_n}) 
            - C(\sigma_n, u_{\sigma_n})} < \omega_n < C(\sigma_n,
            \mu_{\sigma^-_n}) \vee C(\sigma_n, u_{\sigma_n})$,} \\
        \partial   &\text{else.} \\
    \end{dcases}
\end{equation}
Here, $\ell'_n = \partial$ encodes that the offspring is thinned
for both procedures. Note that, with this definition, the law of
$\ell'_n$ is as prescribed in \eqref{eq:lawCoupling}.
Let us also introduce
\[
    \rho_n = \indic \big\{
        \abs{C(\sigma_n, \mu_{\sigma^-_n}) - C(\sigma_n, u_{\sigma_n})} 
        < \omega_n \le \tfrac{L}{N}I^\eta(\sigma^-_n) + \eta
    \big\},
\]
which is the indicator of the event that a $\dag$ particle needs to be immigrated.
We let
\begin{multline} \label{eq:coupling11Bis}
    R_{n+1} = R_n - \delta_{(\sigma_n,\ell_n)} + 
    \indic_{\{ \ell'_n \ne \partial\}} \sum_{a \in \mathcal{P}_n} \delta_{(\sigma_n, \ell'_n)}
    \\
    + \indic_{\{ \ell'_n \ne (1,1)\}} \sum_{a \in \mathcal{P}^*_n}
    \delta_{(\sigma_n+a, *)} 
    + \rho_n \sum_{a \in \mathcal{P}^\dag_n} \delta_{(\sigma_n+a, \dag)}.
\end{multline}
and
\[
    L_{n+1} = L_n + \indic_{\{ \ell'_n \ne \partial\}} \delta_{(\sigma_n, \ell'_n)}
    + \indic_{\{ \ell'_n \ne (1,1)\}} \delta_{(\sigma_n, *)} 
    + \rho_n \delta_{(\sigma_n, \dag)}.
\]

Consider the particles of type $(1,0)$ and $(1,1)$, which are those recorded in
$\mu_{\sigma^-_n}$. No such particle can be born of a $*$, $\dag$, or $(0,1)$
parent, so that we can discard these particles from the dynamics.
Clearly, by \eqref{eq:couplingDag}, the birth at time $t$ from a $(1,0)$
parent yields a $(1,0)$ offspring with probability $C(t, \mu_{t^-})$.
Similarly, if the parent has label $(1,1)$, \eqref{eq:coupling11} ensures
that the offspring has type $\ell' \in \{(1,0), (1,1)\}$ with probability 
$C(t, \mu_{t^-})$. Therefore, irrespective of the type of the parent,
the offspring of a $(1,0)$ or $(1,1)$ parent is of type $(1,0)$ or
$(1,1)$ with probability $C(t, \mu_{t^-})$, which is the construction of
the CMJ process with interaction. By the same argument applied to
$(1,1)$ and $(0,1)$, $\overline{\mu}_t$ is the empirical age distribution
of $N$ i.i.d.\ copies of the non-linear CMJ tree.

Note that, formally, our algorithmic construction only provides us with
the age structures of the CMJ processes, and not with their tree
structures. However, it is straightforward to adapt the construction to
keep track of the Ulam--Harris labels in such a way that the forest
spanned by $(1,0)$ (resp.\ $(0,1)$) and $(1,1)$ particles is distributed
as a forest of $N$ CMJ trees with interaction (resp.\ non-linear CMJ
trees). We chose to not include this to reduce the bookkeeping to a minimum.

The comparison with the process with immigration is more complex. First,
one can directly see from \eqref{eq:coupling11Bis} that, at each potential
birth of a $(1,1)$ parent, either a $(1,1)$ or a $*$ particle is added to
the left of Figure~\ref{fig:coupling}. Moreover, by
\eqref{eq:couplingStar} the birth of a $*$ parent is always effective and
yields a $*$ offspring. Thus the process counting the number of $*$ and
$(1,1)$ particles is an un-thinned CMJ process. The same arguments
applies to $(1,0)$, $(0,1)$, and $\dag$ particles: the tree grown out of
any such particle is distributed as a CMJ tree.

Finally, we need to consider the immigration of $(1,0)$, $(0,1)$, or $\dag$
particles from a $(1,1)$ or $*$ parent. Let us introduce a stopping time 
\[
    \tau = \inf \{ n \ge 0 :
        \abs{C(\sigma_n, \mu_{\sigma^-_n}) - C(\sigma_n, u_{\sigma_n})} >
        \tfrac{L}{N}I^\eta(\sigma^-_n) + \eta
    \}.
\]
By \eqref{eq:couplingStar}, a $\dag$ particle is added at the birth from a
$*$ parent at time $t$ with probability $((LI^\eta(t-))/N + \eta) \wedge 1$. If
the parent is of type $\ell = (1,1)$, \eqref{eq:coupling11} shows that
a new particle of type $(1,0)$ or $(0,1)$ is added with probability
$\abs{C(t, \mu_{t^-}) - C(t, u_t)}$. If not such offspring is added, a
new $\dag$ particle is added with probability 
\[
    \big( 1 \wedge \big(\tfrac{LI^\eta(t-)}{N} + \eta \big) - 
    \abs{C(t, \mu_{t^-}) - C(t, u_t)} \big) \vee 0.
\]
We see that a particle of type $(1,0)$, $(0,1)$, or $\dag$ is added at
the birth of a $(1,1)$ parent with overall probability 
\[
    \abs{C(t, \mu_{t^-}) - C(t, u_t)} + 
    \big( 1 \wedge \big(\tfrac{LI^\eta(t-)}{N} + \eta \big) - 
    \abs{C(t, \mu_{t^-}) - C(t, u_t)} \big) \vee 0.
\]
Thus, $I^\eta(t)$ is \emph{not} distributed as a CMJ process with
immigration. However, as long as $t < \sigma_\tau$, this probability
becomes $1 \wedge ((LI^\eta(t^-))/N + \eta)$ and we do recover the transition
of a CMJ process with immigration. From this observation, one can
construct a (true) CMJ process with immigration
$(\widetilde{I}^\eta(t))_{t \ge 0}$ such that 
\[
    \forall t < \sigma_{\tau}, \quad \widetilde{I}^\eta(t) = I^\eta(t).
\]

Let us finally show that \eqref{eq:goodCoupling} holds on 
the event that $\{ L d_{\Pr}(\overline{\mu}^N_T, u_T) < \eta\}$. 
On this event, using the Lipschitz condition \eqref{eq:Lipschitz}, 
for $\sigma_n < T$,
\begin{align*}
    \abs{C(\sigma_n, \mu_{\sigma_n^-}) - C(\sigma_n, u_{\sigma_n})}
    &\le L d_\Pr\big( \mu_{\sigma_n^-}, u_{\sigma_n} \big)
    \le L d_\Pr\big( \mu_{\sigma_n^-}, \overline{\mu}_{\sigma_n^-} \big) 
    + L d_\Pr\big( \overline{\mu}_{\sigma_n^-}, u_{\sigma_n} \big) \\
    &\le \tfrac{L}{N} I^\eta(\sigma_n^-) + \eta
\end{align*}
and thus $\tau > T$. The coupling is constructed in such a way that 
\begin{equation*} 
    \sum_{i=1}^N \sum_{u \in \mathcal{U}_i} 
    \indic_{\{\sigma^N_u \le t, \sigma'_u \le t\}} 
    \indic_{\{ \sigma^N_u \ne \sigma'_u\}} 
    \le I^\eta(t),
\end{equation*}
and $I^\eta(t)$ coincides with a CMJ process with immigration until time
$\sigma_\tau$, hence the result.
\end{proof}

\section{Completing the proofs}
\label{S:proofs}

\subsection{Proof of the main theorem}

\begin{proof}[Proof of Theorem~\ref{thm:main}]
    The result will follow by the law of large numbers and our coupling.
    Let us define
    \[
        \nu_N = \frac{1}{N} \sum_{i=1}^N \delta_{\mathcal{T}_i},
        \qquad
        \nu'_N = \frac{1}{N} \sum_{i=1}^N \delta_{\mathcal{T}'_i},
    \]
    where $\mathcal{T}'_i$ is the non-linear branching process
    grown out of the $i$-th ancestor in the coupling of
    Proposition~\ref{prop:coupling}.

    Recall from \eqref{eq:defTrees} that $\mathscr{T}$ stands for the set
    of trees. By \cite[Lemma~4.1]{kallenberg2017random}, our result is
    proved if we can show that, for any $\phi \colon \mathscr{T} \to \R$
    which is bounded and continuous in the product topology, 
    \begin{equation} \label{eq:mainProof2}
        \angle{\nu_N, \phi} \longrightarrow
        \E[\phi(\mathcal{T}')],
        \qquad \text{as $N \to \infty$,}
    \end{equation}
    in probability. Let us assume first that there exists $K \ge 1$ such
    that 
    \[
        (s_u = s'_u,\, \abs{u} \le K ) \implies 
        \phi\big( (s_u)_{u \in \mathcal{U}} \big)
        =
        \phi\big( (s'_u)_{u \in \mathcal{U}} \big).
    \]
    That is, that $\phi$ only depends on the tree up to height $K$.
    Let us denote by 
    \[
        M_K \colon (s_u)_{u \in \mathcal{U}} \mapsto
        \sup_{\abs{u} \le K} (\indic_{\{ s_u < \infty\}} s_u)
    \]
    the map giving the height of the first $K$ generations of the tree.

    Fix some $T, \eta > 0$. By the triangular inequality,
    \begin{equation*} 
        \abs{\angle{\nu_N, \phi} - \E[\phi(\mathcal{T}')]} \le 
        \abs{\angle{\nu_N', \phi} - \E[\phi(\mathcal{T}')]}
        + \abs{\angle{\nu_N - \nu'_N, \phi}}.
    \end{equation*}
    The first term vanishes almost surely by the law of large numbers.
    We need to control the second terms. 
    Let us denote by $\mathcal{T}'_{i,K}$ and $\mathcal{T}^N_{i,K}$ the
    restrictions of $\mathcal{T}'_i$ and $\mathcal{T}^N_i$ to the first
    $K$ generations, respectively. To obtain an upper bound for the
    second term, we first use that
    \[
        \abs{\angle{\nu_N - \nu'_N, \phi}}
        \le 
        \frac{1}{N} \norm{\phi}_\infty 
        \Card \big\{
            i \le N : \mathcal{T}'_{i,K} \ne \mathcal{T}^N_{i,K}
        \big\}.
    \]
    Moreover, by a union bound,
    \begin{multline} \label{eq:mainProof1}
        \Card \big\{
            i \le N : \mathcal{T}'_{i,K} \ne \mathcal{T}^N_{i,K}
        \big\} \\
        \le 
        \Card \big\{
            i \le N : M_K(\mathcal{T}'_i) > T
        \big\} 
        + \Card \big\{
            i \le N : M_K(\mathcal{T}^N_i) > T
        \big\} \\
        + \Card \big\{
            i \le N : M_K(\mathcal{T}'_i) \vee M_K(\mathcal{T}^N_i) \le T, 
            \mathcal{T}'_{i,K} \ne \mathcal{T}^N_{i,K}
        \big\}.
    \end{multline}

    For the first two terms, note that both $\mathcal{T}'_i$ and
    $\mathcal{T}^N_i$ are dominated by the un-thinned CMJ tree
    $\mathcal{T}_i$. Therefore, by the law of large numbers,
    \begin{multline*}
        \frac{1}{N} \Card \big\{
            i \le N : M_K(\mathcal{T}'_i) > T
        \big\} 
        + 
        \frac{1}{N} \Card \big\{
            i \le N : M_K(\mathcal{T}^N_i) > T
        \big\} \\
        \le \frac{2}{N} \Card \big\{
            i \le N : M_K(\mathcal{T}_i) > T
        \big\}
        \longrightarrow
        2\P( M_K(\mathcal{T}_1) > T), \qquad \text{as $N \to \infty$}.
    \end{multline*}
    This term will be made arbitrarily small by letting $T \to \infty$.

    For the remaining term in \eqref{eq:mainProof1}, note that    \begin{equation*}
        \Card \big\{
            i \le N : M_K(\mathcal{T}'_i) \vee M_K(\mathcal{T}^N_i) \le T, 
            \mathcal{T}'_{i,K} \ne \mathcal{T}^N_{i,K}
        \big\} 
        \le 
        \sum_{i=1}^N \sum_{u \in \mathcal{U}_i} 
        \indic_{\{ \sigma^N_u \le T, \sigma'_u \le T \}} 
        \indic_{\{ \sigma^N_u \ne \sigma'_u \}}.
    \end{equation*}

    Fix $\epsilon > 0$. Putting the previous two steps together, by
    Proposition~\ref{prop:coupling}, for any $T > 0$ and $\eta > 0$
    \begin{multline*}
        \P\Big( \frac{1}{N} \Card \big\{ i \le N : \mathcal{T}'_{i,K} \ne
            \mathcal{T}^N_{i,K} \big\} \ge 2\epsilon \Big)
        \le 
        \P\Big( \frac{2}{N} \Card\big\{ i \le N : M_K(\mathcal{T}_i) > T
        \big\} \ge \epsilon \Big)
        \\
        + 
        \P\Big( \frac{I^\eta(T)}{N} \ge \epsilon \Big) + \P\big(
        d_\Pr(\overline{\mu}^N_T, u_T) \ge \eta \big).
    \end{multline*}
    By choosing $T$ such that $2\P(M_K(\mathcal{T}_1) > T) < \epsilon$,
    the law of large numbers entails that 
    \[
        \limsup_{N \to \infty}
        \P\Big( \frac{1}{N} \Card \big\{ i \le N : \mathcal{T}'_{i,K} \ne
            \mathcal{T}^N_{i,K} \big\} \ge 2\epsilon \Big)
        \le 
        \limsup_{N \to \infty} \P( I^\eta(T) > \epsilon N),
    \]
    and the claim follows by letting $\eta \to 0$ using
    Lemma~\ref{lem:convImmigration}.

    Up to this point we have shown that, in probability for the weak
    topology,
    \begin{equation} \label{eq:mainProof3}
        \forall K \ge 1,\quad (f_K)_* \nu_N \longrightarrow
        \mathscr{L}(f_K(\mathcal{T}')),\qquad \text{as $N \to \infty$},
    \end{equation}
    where $f_K \colon (s_u)_{u \in \mathcal{U}} \mapsto (s_u)_{\abs{u} \le
    K}$ cuts a tree at height $K$ and $(f_K)_*\nu$ is the corresponding
    pushforward measure.
    Since a measure on $(\R_+ \cup \{ \infty \})^{\mathcal{U}}$
    converges weakly in the product topology if and only if its
    finite-dimensional marginals converge weakly,
    \eqref{eq:mainProof3} directly entails that 
    \[
        \nu_N \longrightarrow \mathscr{L}(\mathcal{T}'), 
        \qquad \text{ as $N \to \infty$,}
    \]
    in probability. 
\end{proof}

\subsection{Limit of the age structure}

\begin{proof}[Proof of Corollary~\ref{cor:ageStructure}]
    Fix a continuous, bounded $\phi \colon \R_+ \to \R$, and let
    \[
        F^\phi_t \colon (s_u)_{u \in \mathcal{U}}
        \mapsto
        \sum_{u \in \mathcal{U}} \indic_{\{ s_u \le t \}} \phi(t - s_u).
    \]
    With this definition,
    \[
        \angle{\mu^N_t, \phi} = \frac{1}{N} \sum_{i=1}^N \sum_{u \in
        \mathcal{U}_i} \indic_{\{ \sigma^N_u \le t \}} \phi(t - \sigma^N_u)
        = \angle{\nu^N, F^\phi_t}.
    \]
    Since $F^\phi_t$ is continuous for almost every realisation of
    $\mathcal{T}'$, Theorem~\ref{thm:main} shows that
    \[
        \angle{\mu_t^N, \phi} = \angle{\nu^N, F^\phi_t}
        \longrightarrow \E[ F^\phi_t(\mathcal{T}') ]
        = 
        \E[ \angle{\mu'_t, \phi} ] = \angle{u_t, \phi}
    \]
    in probability as $N \to \infty$.
    Therefore, by another application of \cite[Lemma~4.1]{kallenberg2017random},
    $\mu^N_t(\diff a) \to u_t(a) \diff a$ in probability for the weak
    topology, for all $t \ge 0$. Since the limit is deterministic, this
    shows finite-dimensional convergence of the process $(\mu^N_t)_{t \ge
    0}$. It remains to prove tightness.

    By \cite[Theorem~16.27]{Kallenberg2002}, it is sufficient to 
    prove tightness of the real-valued process 
    \[
        \forall t \ge 0,\quad \mu^\phi_t = \angle{\mu^N_t, \phi},
    \]
    for each uniformly continuous, bounded $\phi \colon \R_+ \to \R$. For
    $s \le t$ such that $t-s \le \delta$,
    \begin{align*}
        \abs{\mu^\phi_t - \mu^\phi_s}
        &\le \frac{1}{N} \sum_{i=1}^N 
        \sum_{u\in\mathcal{U}_i} 
        \abs*{\indic_{\{ \sigma^N_u \le t\}} \phi(t-\sigma^N_u)
            - \indic_{\{ \sigma^N_u \le s\}} \phi(s-\sigma^N_u)} \\
        &\le 
        \begin{multlined}[t]
        \frac{1}{N} \sum_{i=1}^N \sum_{u \in \mathcal{U}_i}
        \indic_{\{ \sigma^N_u \le s \}} 
        \abs{ \phi(t-\sigma^N_u) - \phi(s-\sigma^N_u)} \\
        + \frac{\norm{\phi}_\infty}{N} \sum_{i=1}^N \sum_{u \in \mathcal{U}_i} 
        \indic_{\{ s < \sigma^N_u \le t\}}
        \end{multlined} \\
        &\le \frac{w_\delta(\phi)}{N} \sum_{i=1}^N Z_i(t) 
        + \frac{\norm{\phi}_\infty}{N} \sum_{i=1}^N (Z_i(t)-Z_i(s)),
    \end{align*}
    where $w_\delta(\phi)$ is the modulus of continuity of $\phi$
    and 
    \[
        \forall t \ge 0,\quad Z_i(t) = \sum_{u \in \mathcal{U}_i}
        \indic_{\{ \sigma_u \le t \}}
    \]
    is the number of individuals at time $t$ of the (un-thinned) CMJ
    process grown out of the $i$-th ancestor. If $\bar{Z}^N_t = (Z_1(t) +
    \dots + Z_n(t)) / N$, for any $T > 0$,
    \begin{equation} \label{eq:proofCor1}
        \E\big[ w_\delta\big((\mu^\phi_t)_{t \le T}\big) \wedge 1 \big] \le 
        w_\delta(\phi) \E[Z_1(T)] + \norm{\phi}_\infty
        \E\big[ w_\delta\big( (\bar{Z}^N_t)_{t \le T} \big) \wedge 1 \big].
    \end{equation}
    By the law of large numbers, $\bar{Z}^N(t) \to \E[Z_1(t)]$ almost
    surely. Moreover, since this is an increasing process and $t \mapsto
    \E[Z_1(t)]$ is continuous, \cite[Theorem 3.37, Chapter VI]{jacod2013limit}
    shows that $(\bar{Z}^N_t)_{t \ge 0}$ converges to its limit in the
    topology of uniform convergence on compact sets. Thus, 
    by \cite[Theorem~16.5]{Kallenberg2002} applied to $(\bar{Z}^N_t)_{t
    \le T}$ and \eqref{eq:proofCor1}, 
    \[
        \adjustlimits
        \lim_{\delta \to 0} \limsup_{N \to \infty} 
        \E\big[ w_\delta\big( (\mu_t^\phi)_{t \le T} \big) \wedge 1 \big] 
        = 0,
    \]
    which (again by \cite[Theorem~16.5]{Kallenberg2002}) entails that
    $(\mu^\phi_t)_{t \ge 0}$ is tight in the topology of uniform
    convergence on compact sets.
\end{proof}

\subsection{Limit of the infection chains}

As indicated in the introduction, we prove Corollary~\ref{cor:infectionChain} 
by writing an infection chain as a functional of a tree, then use the
continuous mapping theorem, and evaluate the functional for the
non-linear branching process. We start by computing the functional at the
limit.

Recall the Markov chain $(T_i)_i$ defined above Corollary~\ref{cor:infectionChain}.
We think of it as a random element of $\cup_{k \ge 1} \R^k$, and write
$\mathbf{t} = (t_i)_{i=1}^k$ for a generic element of this set. Recall
that, for $u \in \mathcal{U}_1$, we write $(\sigma_v)_{1 \preceq v \preceq u}$ 
for the random collection of birth times leading to $v$, which has length
$\abs{u}$ and is indexed backward-in-time from $u$ to the ancestor $1$.

\begin{lemma} \label{lem:ancestralNonLinear}
    Let $\mathcal{T}' = (\sigma'_u)_{u \in \mathcal{U}_1}$ be the
    non-linear CMJ tree. Then
    \[
        \forall T > 0,\quad 
        \E\Big[ \sum_{u \in \mathcal{U}_1} \indic_{\{ \sigma'_u \le T\}} 
        \delta_{(\sigma'_v)_{1 \preceq v \preceq u}}(\diff \mathbf{t}) \Big]
        = 
        \int_{-\infty}^T u_s(0) \P\big( (T_i)_i \in \diff \mathbf{t} \mid T_1 = s \big)
        \diff s.
    \]
\end{lemma}

\begin{proof}
For $\mathbf{t} = (t_i)_{i=1}^k$, with $k \ge 1$ and  $t_1 \ge \dots \ge
t_{k-1} \ge 0 > t_k$, let $P_k(\mathbf{t})$ be the mean density of
individuals with infection chain $\mathbf{t}$ with length $k$. Formally,
for any non-negative test function $\phi \colon \R^k \to \R_+$,
$P_k(\mathbf{t})$ is defined as
\[
    \int_{\R_-} \Big( \int_{\R_+} P_k\big( \mathbf{t} \big) 
        \phi\big( \mathbf{t} \big) \diff t_1 \dots \diff t_{k-1}
    \Big) \diff t_k
    =
    \E\Big[ 
        \sum_{u \in \mathcal{U}_1, \abs{u}=k} 
        \indic_{\{ \sigma'_u < \infty \}}
        \phi\big( (\sigma'_v)_{1 \preceq v \preceq u} \big) 
    \Big].
\]
(The existence of such a density follows directly from the existence of a
density $\tau$ for the mean measure $\E[\mathcal{P}(\diff a)]$.)
Note that, with this definition, $P_1(t_1)$ is the density of infection
times of the initial individuals and thus $P_1(t_1) = u_0(-t_1)$. Using
the branching property, we can compute recursively that
\begin{align*}
    \E\Big[ 
        \sum_{\abs{u}=k+1} 
        &\indic_{\{ \sigma'_u < \infty \}}
        \phi\big( (\sigma'_v)_{1 \preceq v \preceq u} \big) 
    \Big] \\
    &=
    \E\Big[ 
        \sum_{\abs{u}=k} 
        \sum_{j \ge 1}
        \indic_{\{ \sigma'_{uj} < \infty \}}
        \phi\big( \sigma'_{uj}, (\sigma'_v)_{1 \preceq v \preceq u} \big) 
    \Big] \\
    &=
    \E\Big[ 
        \sum_{\abs{u}=k} 
        \indic_{\{ \sigma'_u < \infty \}}
        \int_0^\infty
        \phi\big( \sigma'_u+a, (\sigma'_v)_{1 \preceq v \preceq u} \big)  
        \tau(a) C(\sigma'_u+a, u_{\sigma'_u+a}) \diff a
    \Big].
\end{align*}
This yields the recursion 
\[
    P_{k+1}\big( (t_i)_{i=1}^{k+1} \big) = P_k\big( (t_{i+1})_{i=1}^k
    \big) \tau(t_1-t_2)
    C\big(t_1, u_{t_1}\big).
\]
Hence, using that $u_t(s) = u_{t-s}(0)$,
\begin{align*}
    P_k\big( (t_i)_{i=1}^k \big) 
    &= P_1(t_k) \prod_{i=1}^{k-1} C(t_i, u_{t_i}) \tau(t_i-t_{i+1}) \\
    &= u_0(-t_k) \prod_{i=1}^{k-1} C(t_i, u_{t_i}) \tau(t_i-t_{i+1})
    \frac{u_{t_i}(t_i-t_{i+1})}{u_{t_{i+1}}(0)} \\
    &= \frac{u_{t_1}(0)}{u_{t_k}(0)} 
    u_0(-t_k) \prod_{i=1}^{k-1} \frac{C(t_i, u_{t_i})}{u_{t_i}(0)}
    u_{t_i}(t_i-t_{i+1}) \tau(t_i-t_{i+1}) \\
    &= u_{t_1}(0) 
    \P\big( (T_i)_{i=1}^k = (t_i)_{i=1}^k \mid T_1 = t_1 \big),
\end{align*}
where from the second line on we have used the convention that $u_{-s}(0) =
u_0(s)$ for $s \ge 0$. The result is obtained by summing over $k$ and
integrating over $t_1 \le T$.
\end{proof}

\begin{proof}[Proof of Corollary~\ref{cor:infectionChain}]
    Let $\phi \colon \cup_{k \ge 1} \R^k \to \R_+$ be a non-negative
    continuous bounded function and, for $u \in \mathcal{U}$, let $\phi_u$
    be the functional
    \[
        \phi_u \colon (s_v)_{v \in \mathcal{U}} \mapsto 
        \phi\big( (s_v)_{\varnothing \preceq v \preceq u} \big)
    \]
    which is clearly also continuous bounded. Recall that viewing a tree
    as an element of $\mathscr{T}$ required us to relabel the individuals
    so that the root is $\varnothing$. Therefore, infection chains in
    $\mathscr{T}$ are naturally indexed from $u$ to $\varnothing$,
    whereas they are indexed from $u$ to $i$ in $\mathcal{T}_i$.

    By the continuous mapping theorem and Theorem~\ref{thm:main},
    \[
        \frac{1}{N} \sum_{i=1}^N \phi\big( 
            (\sigma^N_{(i,v)})_{\varnothing \preceq v \preceq u} 
        \big)
        \indic_{\{ \sigma^N_{(i,u)} \le T \}}
        \longrightarrow
        \E\big[ \indic_{\{ \sigma'_{(1,u)} \le T\}} 
        \phi\big( (\sigma'_{(1,v)})_{\varnothing \preceq v \preceq u} \big) \big]
    \]
    in distribution as $N \to \infty$. By Skorohod's representation
    theorem (for instance \cite[Theorem~6.7]{billingsley1999}), let us
    assume that this convergence holds almost surely for all $i \le N$
    and $u \in \mathcal{U}$. Let us show that 
    \begin{equation} \label{eq:proofCor2ToShow}
        \frac{1}{N} \sum_{i=1}^N \sum_{u \in \mathcal{U}_i}
        \indic_{\{ \sigma^N_u \le T \}} \phi\big( (\sigma^N_v)_{i \preceq v \preceq u} \big)
        \longrightarrow
        \sum_{u \in \mathcal{U}_1} 
        \E\big[ \indic_{\{\sigma'_u \le T\}}
            \phi\big( (\sigma'_v)_{1 \preceq v \preceq u} \big)
        \big], \qquad \text{as $N \to \infty$,}
    \end{equation}
    almost surely, by applying the dominated convergence theorem to each
    realisation. By dominating the CMJ tree with interaction by that
    without thinning, for any $u \in \mathcal{U}$,
    \begin{equation} \label{eq:proofCor2Bound}
        \frac{1}{N} \sum_{i=1}^N \indic_{\{ \sigma^N_{(i,u)} \le T \}} 
        \phi\big( (\sigma^N_{(i,v)})_{\varnothing \preceq v \preceq u} \big)
        \le 
        \frac{\norm{\phi}_\infty}{N} 
        \sum_{i=1}^N \indic_{\{\sigma_{(i,u)} \le T\}}.
    \end{equation}
    By the law of large numbers, almost surely as $N \to \infty$,
    \[
        X_u^N \coloneqq \frac{1}{N} \sum_{i=1}^N \indic_{\{\sigma_{(i,u)} \le T\}}
        \longrightarrow
        \E[ \indic_{\{\sigma_{(1,u)} \le T\}} ]
    \]
    and 
    \[
        \sum_{u \in \mathcal{U}} X_u^N =
        \frac{1}{N} \sum_{i=1}^N \sum_{u \in \mathcal{U}_i}
        \indic_{\{\sigma_u \le T\}}
        \longrightarrow
        \E[ Z_1(t) ] = 
        \sum_{u \in \mathcal{U}} 
        \E\big[ \indic_{\{\sigma_{(1,u)} \le T\}} \big].
    \]
    These two points show that, for a.e.\ realisation, the sequence 
    $( (X_u^N)_{u \in \mathcal{U}} )_{N \ge 1}$ is uniformly integrable
    as a sequence of maps $X^N \colon u \mapsto X^N_u$. (This follows for
    instance from a simple adaptation of \cite[Lemma~4.1]{kallenberg2017random}.) 
    Therefore, the variables on the left-hand side of \eqref{eq:proofCor2Bound}
    are also a.s.\ uniformly integrable (as functions of $u \in
    \mathcal{U}$), and we can apply the bounded convergence theorem to
    a.e.\ realisation to deduce \eqref{eq:proofCor2ToShow}. Finally,
    Lemma~\ref{lem:ancestralNonLinear} proves that the right-hand side of
    \eqref{eq:proofCor2ToShow} coincides with the limit in the statement of
    Corollary~\ref{cor:infectionChain}.
\end{proof}

\subsection*{Acknowledgments}

The authors would like to express their gratitude to Jean-Jil Duchamps,
who was involved in several discussions at an early stage of this work.

\bibliographystyle{plain}
\bibliography{contacts}

\begin{thebibliography}{10}

\bibitem{Barbour2013}
Andrew Barbour and Gesine Reinert.
\newblock Approximating the epidemic curve.
\newblock {\em Electronic Journal of Probability}, 18(none):1--30, Jan 2013.

\bibitem{billingsley1999}
Patrick Billingsley.
\newblock {\em Convergence of probability measures}.
\newblock Wiley Series in Probability and Statistics: Probability and Statistics. John Wiley \& Sons Inc., New York, second edition, 1999.
\newblock A Wiley-Interscience Publication.

\bibitem{britton2019estimation}
Tom Britton and Gianpaolo Scalia~Tomba.
\newblock Estimation in emerging epidemics: biases and remedies.
\newblock {\em Journal of the Royal Society Interface}, 16(150):20180670, 2019.

\bibitem{calvez2022dynamics}
Vincent Calvez, Beno{\^\i}t Henry, Sylvie M{\'e}l{\'e}ard, and Viet~Chi Tran.
\newblock Dynamics of lineages in adaptation to a gradual environmental change.
\newblock {\em Annales Henri Lebesgue}, 5:729--777, 2022.

\bibitem{Chaintron2022}
Louis-Pierre Chaintron and Antoine Diez.
\newblock Propagation of chaos: A review of models, methods and applications. i. models and methods.
\newblock {\em Kinetic and Related Models}, 15(6):895--1015, 2022.

\bibitem{chevallier2017}
Julien Chevallier.
\newblock Mean-field limit of generalized hawkes processes.
\newblock {\em Stochastic Processes and their Applications}, 127(12):3870--3912, 2017.

\bibitem{Coste2021}
Christophe F.~D. Coste, Fran{\c{c}}ois Bienvenu, Victor Ronget, Juan-Pablo Ramirez-Loza, Sarah Cubaynes, and Samuel Pavard.
\newblock The kinship matrix: inferring the kinship structure of a population from its demography.
\newblock {\em Ecology Letters}, 24(12):2750--2762, 2021.

\bibitem{dewitt2024mean}
William~S DeWitt, Steven~N Evans, Ella Hiesmayr, and Sebastian Hummel.
\newblock Mean-field interacting multi-type birth--death processes with a view to applications in phylodynamics.
\newblock {\em Theoretical Population Biology}, 159:1--12, 2024.

\bibitem{duchamps2023general}
Jean-Jil Duchamps, F{\'e}lix Foutel-Rodier, and Emmanuel Schertzer.
\newblock General epidemiological models: Law of large numbers and contact tracing.
\newblock {\em Electronic Journal of Probability}, 28:1--37, 2023.

\bibitem{favero2022modelling}
Martina Favero, Gianpaolo Scalia~Tomba, and Tom Britton.
\newblock Modelling preventive measures and their effect on generation times in emerging epidemics.
\newblock {\em Journal of The Royal Society Interface}, 19(191):20220128, 2022.

\bibitem{forien2025stochastic}
Rapha{\"el} Forien, Guodong Pang, {\'Etienne}~Pardoux, and Ars{\`ene}~Brice Zotsa-Ngoufack.
\newblock Stochastic epidemic models with varying infectivity and waning immunity.
\newblock {\em arXiv}, 2022.

\bibitem{foutel2022individual}
F{\'e}lix Foutel-Rodier, Fran{\c{c}}ois Blanquart, Philibert Courau, Peter Czuppon, Jean-Jil Duchamps, Jasmine Gamblin, {\'E}lise Kerdoncuff, Rob Kulathinal, L{\'e}o R{\'e}gnier, Laura Vuduc, et~al.
\newblock From individual-based epidemic models to mckendrick-von foerster pdes: A guide to modeling and inferring covid-19 dynamics.
\newblock {\em Journal of Mathematical Biology}, 85(4):43, 2022.

\bibitem{henry2023time}
Beno{\^\i}t Henry, Sylvie M{\'e}l{\'e}ard, and Viet~Chi Tran.
\newblock Time reversal of spinal processes for linear and non-linear branching processes near stationarity.
\newblock {\em Electronic Journal of Probability}, 28:1--27, 2023.

\bibitem{jacod2013limit}
Jean Jacod and Albert Shiryaev.
\newblock {\em Limit theorems for stochastic processes}, volume 288.
\newblock Springer Science \& Business Media, 2013.

\bibitem{jagers1975branching}
Peter Jagers.
\newblock {\em Branching processes with biological applications}.
\newblock Wiley-Interscience [John Wiley \& Sons], London, 1975.

\bibitem{jagers1984growth}
Peter Jagers and Olle Nerman.
\newblock The growth and composition of branching populations.
\newblock {\em Advances in applied probability}, 16(2):221--259, 1984.

\bibitem{Kallenberg2002}
Olav Kallenberg.
\newblock {\em Foundations of Modern Probability}.
\newblock Springer New York, New York, NY, 2002.

\bibitem{kallenberg2017random}
Olav Kallenberg et~al.
\newblock {\em Random measures, theory and applications}, volume~1.
\newblock Springer, 2017.

\bibitem{lyons1995conceptual}
Russell Lyons, Robin Pemantle, and Yuval Peres.
\newblock Conceptual proofs of {$L \log L$} criteria for mean behavior of branching processes.
\newblock {\em The Annals of Probability}, pages 1125--1138, 1995.

\bibitem{marguet2019uniform}
Aline Marguet.
\newblock Uniform sampling in a structured branching population.
\newblock {\em Bernoulli}, 25(4A):2649--2695, 2019.

\bibitem{Nerman1981}
Olle Nerman.
\newblock On the convergence of supercritical general ({C-M-J}) branching processes.
\newblock {\em Zeitschrift f{\"u}r Wahrscheinlichkeitstheorie und Verwandte Gebiete}, 57(3):365--395, 1981.

\bibitem{Nerman1984}
Olle Nerman and Peter Jagers.
\newblock The stable doubly infinite pedigree process of supercritical branching populations.
\newblock {\em Zeitschrift f{\"u}r Wahrscheinlichkeitstheorie und Verwandte Gebiete}, 65(3):445--460, 1984.

\bibitem{overbeck1996nonlinear}
L~Overbeck.
\newblock Nonlinear superprocesses.
\newblock {\em The Annals of Probability}, 24(2):743--760, 1996.

\bibitem{pakkanen2023unifying}
Mikko~S Pakkanen, Xenia Miscouridou, Matthew~J Penn, Charles Whittaker, Tresnia Berah, Swapnil Mishra, Thomas~A Mellan, and Samir Bhatt.
\newblock Unifying incidence and prevalence under a time-varying general branching process.
\newblock {\em Journal of Mathematical Biology}, 87(2):35, 2023.

\bibitem{Shi2015}
Zhan Shi.
\newblock {\em Branching Random Walks: {\'E}cole d'{\'E}t{\'e} de Probabilit{\'e}s de Saint-Flour XLII -- 2012}.
\newblock Springer International Publishing, Cham, 2015.

\bibitem{sznitman1991}
Alain-Sol Sznitman.
\newblock Topics in propagation of chaos.
\newblock In Paul-Louis Hennequin, editor, {\em Ecole d'Et{\'e} de Probabilit{\'e}s de Saint-Flour XIX --- 1989}, pages 165--251, Berlin, Heidelberg, 1991. Springer Berlin Heidelberg.

\bibitem{webb1985theory}
Glenn~F Webb.
\newblock {\em Theory of nonlinear age-dependent population dynamics}.
\newblock CRC Press, 1985.

\end{thebibliography}

\end{document}